\documentclass[aap]{imsart}

\RequirePackage{amsthm,amsmath,amsfonts,amssymb}
\RequirePackage[numbers]{natbib}
\usepackage{amsfonts}
\usepackage{amsmath}
\usepackage{amsthm}
\usepackage{bm,bbm}
\usepackage[dvips]{graphicx}
\usepackage{caption}
\usepackage{color}

\usepackage{bigints}

\usepackage{mathtools}
\usepackage{commath}

\startlocaldefs
\theoremstyle{plain}
\newtheorem{theorem}{Theorem}[section]
\newtheorem{lemma}[theorem]{Lemma}
\newtheorem{proposition}[theorem]{Proposition}
\newtheorem{corollary}[theorem]{Corollary}

\theoremstyle{remark}

\newtheorem{remark}[theorem]{Remark}
\newtheorem{example}[theorem]{Example}

\numberwithin{equation}{section}

\allowdisplaybreaks

\def\f{\frac}
\def\ms{\mathsf}
\def\beqn{\begin{equation}}
\def\beqn*{$$}
\def\eeqn{\end{equation}}

\newcommand{\bx}{{\bf x}}
\newcommand{\by}{{\bf y}}
\newcommand{\bz}{{\bf z}}

\def\P{\mathbb{P}}

\def\E{\mathbb{E}}
\def\Pn{\mathcal P_n}
\def\Bn{\mathcal B_n}

\newcommand{\reals}{{\mathbb R}}
\newcommand{\bbr}{\reals}
\newcommand{\R}{\reals}
\newcommand{\bbn}{{\mathbb N}}

\newcommand{\vep}{\varepsilon}
\newcommand{\bbz}{\protect{\mathbb Z}}

\newcommand{\X}{{\mathcal{X}}}
\newcommand{\Y}{{\mathcal{Y}}}
\newcommand{\I}{\mathcal I}

\newcommand{\B}{\mathcal B}

\newcommand{\bell}{{\bm \ell}}

\newcommand{\one}{{\mathbbm 1}}

\newcommand{\remove}[1]{}

\newcommand{\C}{\check{C}}

\newcommand{\Z}{\mathcal Z}

\def\enp{\end{proof}}
\def\bel{\begin{lemma}}
\def\bep{\begin{proof}}
\def\enl{\end{lemma}}
\newcommand{\M}{\mathcal M}

\newcommand{\Pnr}{\mathcal P_n|}

\newcommand{\bQ}{Q^\partial}
\newcommand{\bt}{{\bf t}}

\newcommand{\ba}{{\bf a}}

\newcommand{\U}{\mathcal U}
\newcommand{\hi}{h^{(i)}}

\newcommand{\vni}{v_n^{(i)}}
\newcommand{\uni}{u_n^{(i)}}
\newcommand{\bfz}{{\bf 0}}
\newcommand{\weta}{\widetilde \eta}
\newcommand{\zetal}{\zeta^{(\ell)}}
\newcommand{\heta}{\hat \eta}
\newcommand{\hzeta}{\hat \zeta}
\newcommand{\hzetal}{\hat \zeta^{(\ell)}}

\newcommand{\diam}{\mathsf{diam}}

\endlocaldefs

\begin{document}

\begin{frontmatter}

\title{Large deviation principle for geometric and topological functionals and  associated point processes}
	\runtitle{LDP for geometric and topological functionals and point processes}

\begin{aug}
\author[A]{\fnms{Christian} \snm{Hirsch}\ead[label=e1]{hirsch@math.au.dk}},
\and
\author[B]{\fnms{Takashi} \snm{ Owada}\ead[label=e2]{owada@purdue.edu}}
\address[A]{Department of Mathematics\\ Aarhus University \\  Ny Munkegade, 118, 8000, Aarhus C,  Denmark, \printead{e1}}

\address[B]{Department of Statistics\\
Purdue University \\
West Lafayette, 47907, USA, \printead{e2}}
\end{aug}

\begin{abstract}
We prove a large deviation principle for the point process {associated to $k$-element connected components in $\R^d$ with respect to the connectivity radii $r_n\to\infty$}. The random points are generated from a homogeneous Poisson point process {or the corresponding binomial point process}, so that $(r_n)_{n\ge1}$ satisfies $n^kr_n^{d(k-1)}\to\infty$ and $nr_n^d\to0$ as $n\to\infty$ (i.e., sparse regime). The rate function for the obtained large deviation principle can be represented  as relative entropy. As an application, we deduce large deviation principles for various functionals and point processes appearing in stochastic geometry and topology. As concrete examples of topological invariants, we consider persistent Betti numbers of geometric complexes and the number  of  Morse critical points of the min-type distance function. 
\end{abstract}



\end{frontmatter}


\section{Introduction}

The objective of this paper is to examine large deviation behaviors of a point process associated to a configuration of random points generated by a homogeneous Poisson point process {or the corresponding binomial point process}. The \emph{Donsker--Varadhan large deviation principle (LDP)} (see \cite{dembo:zeitouni:1998})  characterizes the limiting behavior of a family of probability measures $(\mu_n)_{n\ge1}$ on a measurable space $(\mathcal X, \B)$. We say that $(\mu_n)_{n\ge1}$ satisfies an LDP with rate $a_n$ and \emph{rate function} $I(x)$, if for any $A\in \B$, 
\begin{equation}  \label{e:intro.LDP}
{-\inf_{x\in \text{int}(A)} I(x)} \le \liminf_{n\to\infty} \frac{1}{a_n}\log \mu_n(A) \le \limsup_{n\to\infty} \frac{1}{a_n}\log \mu_n(A) \le -\inf_{x\in \bar A} I(x), 
\end{equation}
where {$\text{int}(A)$ denotes the interior of $A$} and $\bar A$ the closure of $A$. 

Let $\Pn$ be a homogeneous Poisson point process with intensity $n$ on the unit cube $[0,1]^d$, $d\ge2$, 
and $(r_n)_{n\ge1}$ be a sequence of positive numbers decreasing to $0$ as $n\to\infty$. We put our focus  on the spatial distribution of  $k$-tuples  $\Y\subset \Pn$, which  consist of  components with respect to  $r_n$. More specifically,  given a parameter $t\in (0,\infty)$, we consider the following  point process in $(\R^d)^k$: 
\begin{equation}  \label{e:intro.point.proc}
\sum_{\Y\subset \Pn, \, |\Y|=k} s_n(\Y, \Pn; t)\, \delta_{r_n^{-1} ( \Y-\ell(\Y))}, 
\end{equation}
where $\delta_\bz$ denotes the Dirac measure at $\bz=(z_1,\dots,z_k)\in (\R^d)^k$ and $\ell(\Y)$ is a ``center" point of $\Y$ such as the left most point (in the lexicographic ordering). Moreover, $s_n$ denotes an indicator function,  requiring that the diameter of $k$-tuples $\Y\subset \Pn$ be  at most a constant multiple of $r_n$, and further, such $k$-tuples must {be distant at least} $r_nt$ from all the remaining points in $\Pn$. In this setup, our main theorem (i.e., Theorem \ref{t:LDP.pp.general}) aims to describe  the LDP for the  process \eqref{e:intro.point.proc} in the so-called \emph{sparse regime}: $nr_n^d\to0$ as $n\to\infty$. 

As a primary application, we also provide the LDP for  geometric and topological statistics possessing the structure as  component counts. Given $t\in(0,\infty)$ and measurable functions {$H: (\R^d)^k\to [0, \infty)$} that are symmetric in the $k$ arguments of $\R^d$, we define 
\begin{equation}  \label{e:intro.Gn}
G_n(\Y, \Pn; t) := H(r_n^{-1}\Y)\, \one \big\{ \|y-z\| \ge r_n t, \ \text{for all } y \in \Y \text{ and } z\in \Pn\setminus \Y \big\}, 
\end{equation}
where $\one \{ \cdot \}$ denotes an indicator function and  $\| \cdot \|$ is the Euclidean norm. We then  consider the  point process 
\begin{equation}  \label{e:intro.pp.Ustat}
\sum_{\Y\subset \Pn, \, |\Y|=k} \delta_{G_n(\Y, \Pn; t)}. 
\end{equation}
The process \eqref{e:intro.pp.Ustat} is associated to scaled $k$-tuples $r_n^{-1}\Y$, satisfying the geometric condition implicit in $H$, with the additional restriction  that geometric objects generated by $r_n^{-1}\Y$ {be distant at least} $t$ from all the remaining points of $r_n^{-1}\Pn$. A more precise and general setup for the process \eqref{e:intro.pp.Ustat} is given in Section \ref{sec:application:geom.topo.funcl}. 	

With an appropriate choice of scaling regimes,  the large deviations of point processes as in \eqref{e:intro.pp.Ustat} have been studied  by  Sanov's theorem and its variant \cite{eichelsbacher:lowe:1995, serfling:wang:2000, eichelsbacher:schmock:2002}. In recent times, the  process \eqref{e:intro.pp.Ustat} also found  applications in stochastic geometry \cite{decreusefond:schulte:thaele:2016, otto:2020, bobrowski:schulte:yogeshwaran:2021}. In particular, the authors of  \cite{bobrowski:schulte:yogeshwaran:2021} explored the Poisson process approximation of {the point process for general stabilizing functionals, including \eqref{e:intro.pp.Ustat} as its special case,} by deriving the rate of convergence  in terms of the  Kantorovich-Rubinstein distance. Furthermore, \cite{owada:2022b} discussed  large deviations of the probability distribution of \eqref{e:intro.pp.Ustat}  under the $M_0$-topology, with the assumption that the connectivity radii are even smaller than those considered in \cite{decreusefond:schulte:thaele:2016, bobrowski:schulte:yogeshwaran:2021}. 

{As another main application, we also examine  the LDP for   statistics of the form}
\begin{equation}  \label{e:intro.Ustat}
\sum_{\Y\subset \Pn, \, |\Y|=k} G_n(\Y, \Pn; t).  
\end{equation}
A variety of functionals in stochastic geometry can be considered as such statistics; see \cite{blaszczyszyn:yogeshwaran:yukich:2019,  bobrowski:schulte:yogeshwaran:2021, otto:2020, lachiezerey:reitzner:2016}. 
Moreover, with an appropriate choice of $H$, the statistics \eqref{e:intro.Ustat} can be used to investigate the behavior of topological invariants of a geometric complex \cite{ bobrowski:adler:2014, bobrowski:mukherjee:2015, kahle:meckes:2013,owada:thomas:2020, yogeshwaran:adler:2015, owada:2022a}. Along this line of research, Section \ref{sec:applications} below discusses  the LDP for persistent Betti numbers  and the number of Morse critical points of the min-type distance function. {Loosely speaking, the persistent Betti number is a quantifier for topological complexity, capturing the creation / destruction of topological cycles. The Morse critical points under our consideration provide a good approximation of homological changes in the geometric complex.}

The large deviation behavior of the processes \eqref{e:intro.point.proc}, \eqref{e:intro.pp.Ustat}, and \eqref{e:intro.Ustat} heavily depends on how rapidly the sequence $(r_n)_{n\ge1}$ decays to $0$ as $n\to\infty$. In the existing literature (not necessarily relating to large deviations), the configuration of geometric objects (determined by $H$ in \eqref{e:intro.Gn}) splits into multiple regimes \cite{penrose:2003, kahle:meckes:2013, owada:thomas:2020}. In the context of  random geometric complexes (resp.~random geometric graphs), if $nr_n^d\to0$, called the \emph{sparse regime}, the spatial distribution of complexes (resp.~graphs) is sparse, so that they are mostly observed as isolated components. In the critical phase $nr_n^d\to c\in(0,\infty)$, called the \emph{critical regime} for which $r_n$ decreases to $0$ at a slower rate than the sparse regime, the complex (resp.~graph) begins to coalesce, forming much larger components. Finally, the case when $nr_n^d\to\infty$ is the \emph{dense regime}, for which the complex (resp.~graph) is even more connected and may even consist of a single giant component. 

Among the regimes described above, the present paper focuses on 
 the sparse regime. In this case, there have been a number of studies on the ``average" behavior and the likely deviations from the ``average" behavior of geometric functionals and topological invariants, such as subgraph/component counts and Betti numbers. More specifically, a variety of strong laws of large numbers and central limit theorems for these quantities have already been established. The readers may refer to  the monograph \cite{penrose:2003} for the limit theorems for subgraph/component counts, while  the works in \cite{kahle:meckes:2013, owada:thomas:2020, yogeshwaran:adler:2015} provide the limit theorems  for  topological invariants. 
 
In contrast, however, there have been very few attempts made at examining the  large deviation behaviors of these quantities, especially  from the viewpoints of an LDP. In fact, even for a simple edge count in a random geometric graph, determining the rate $a_n$ and the rate function $I(x)$ in \eqref{e:intro.LDP} is   a highly non-trivial  problem (\cite{chatterjee:harel:2020}). Although there are several works (e.g., \cite{bachmann:reitzner:2018, reitzner:schulte:thale:2017, yogeshwaran:subag:adler:2017}) that deduced concentration inequalities for subgraph counts, the length power functionals, and Betti numbers, these papers were not aimed to derive the LDP. As a consequence,  the obtained upper bounds in these studies do not seem to be tight. Furthermore, \cite{schreiber:yukich:2005} studied the LDP for the functional of spatial point processes satisfying a weak dependence condition {characterized by} a radius of stabilization. One of the major assumptions in their study is that the contribution of any particular vertex must be uniformly bounded (see condition (L1) therein). However,   our functionals and point processes may not fulfill such  conditions.  More importantly, the study in \cite{schreiber:yukich:2005} treated only the critical regime, whereas the main focus of our study is the  sparse regime. One of the benefits of studying the sparse regime is that geometric and topological objects have a relatively simpler structure in the limit, which will allow us to explicitly identify the structure of a rate function. In some cases, we may even solve the variational problem for the rate function; see Remark \ref{rem:explicit.rate.function}. {We also note that \cite{seppalainen:yukich:2001} provides a framework to establish LDPs when the uniform boundedness assumption in \cite{schreiber:yukich:2005} is replaced by the finiteness of suitably scaled logarithmic moment generating functions (see Equ.~(2.5) therein). However, the LDPs in \cite{seppalainen:yukich:2001} are designed only  for real-valued random variables, whereas our main results include LDPs for random measures, such as Theorems \ref{t:LDP.pp.general} and \ref{t:LDP.pp.Ustat} below.} 
As a final remark, we note that the work in \cite{hiraoka:kanazawa:miyanaga:tsunoda:2022+}, which is still under preparation, seems to utilize the ideas of \cite{seppalainen:yukich:2001} to derive the LDP for Betti numbers and persistence diagrams of cubical complexes. 

The remainder of this paper is structured as follows. Section \ref{sec:model.main.results} gives a rigorous description of the LDP for the point process in \eqref{e:intro.point.proc}. In Section \ref{sec:application:geom.topo.funcl}, {we shall deduce}  the LDP for the processes in \eqref{e:intro.pp.Ustat} and \eqref{e:intro.Ustat}. Additional examples on persistent Betti numbers and the number of Morse critical points will be offered in Section \ref{sec:applications}. All the proofs are deferred to Section \ref{sec:proofs}. 
For the proof of our main theorems, we partition the unit cube $[0,1]^d$ into multiple smaller sub-cubes and consider a family of point processes restricted to each of these small sub-cubes. Subsequently, with the aid of the main theorem in \cite{decreusefond:schulte:thaele:2016}, we establish weak convergence of such point processes in terms of the total variation distance. We then exploit Cram\'er's theorem in Polish spaces \cite[Theorem 6.1.3]{dembo:zeitouni:1998} to  identify the structure of a rate function. After that, {Proposition \ref{p:relative.entropy.and.rate}  ensures that the rate function can be represented in terms of relative entropy. 
The required approximation argument relies on the standard  technique on the maximal coupling (\cite[Lemma 4.32]{kallenberg}).} For the proof of the theorems in Section \ref{sec:application:geom.topo.funcl}, we shall utilize an extension of the contraction principle, provided in \cite[Theorem 4.2.23]{dembo:zeitouni:1998}. 

Before concluding the Introduction, we comment on our setup and assumptions. First, we assume that the Poisson point process $\Pn$ is homogeneous. It is clear, however, that in many applications (e.g., \cite{adler:bobrowski:weinberger:2014, owada:adler:2017}),  it is important to understand geometric and topological effects of  lack of homogeneity. A possible starting point for introducing inhomogeneity is to look at point clouds arising from ``inhomogeneous" Poisson point processes. In this case, there is no doubt that a new and more involved machinery must be developed to analyze such data; a detailed discussion will be postponed to a future publication. {Another possible extension seems to investigate the regimes other than the sparse case. At least for the critical case (i.e., $nr_n^d \to c\in (0,\infty)$) however, it is impossible to directly translate our proof techniques. In particular, Lemma \ref{l:total.variation.eta.kn}  does not hold anymore; see Remark \ref{rem:critical.regime} for more details on this point.} 
{Finally},  the LDPs in Section \ref{sec:applications}  are proven only when the dimensions of the Euclidean space, as well as those of topological invariants, are small enough. This is due to the fact that our proof techniques can apply only to a lower-dimensional case. It is still unclear whether the  LDP holds in  higher-dimensions; this is also left as a future topic of our research. 

\section{Model and main results}  \label{sec:model.main.results}

Let $\Pn$ be a homogeneous Poisson point process on $[0,1]^d$, $d\ge2$, with intensity $n$.
{Choose an integer $k \ge 2$, which will remain fixed for the remainder of this section,} and let 
$$
\diam (x_1,\dots,x_k):= \max_{1\le i <j \le k}\|x_i-x_j\|, \ \ \ x_i \in \R^d. 
$$
Taking a sequence $(r_n)_{n\ge1}$ decaying to $0$ as $n\to\infty$, we focus on the \emph{sparse regime}: 
\begin{equation}  \label{e:sparse.radius}
\rho_{k, n}:= n^kr_n^{d(k-1)} \to \infty,   \ \ \ nr_n^d\to0, \ \   \text{ as } n\to\infty. 
\end{equation}
For a subset $\Y$ of $k$ points in $\R^d$, a finite set $\Z \supset \Y$ in $\R^d$ and  {$t\in (0,\infty)$}, we define 
$$
c(\Y, \Z; t) := \one \big\{ \|y-z\| \ge t, \ \text{for all } y \in \Y \text{ and } z\in \Z\setminus \Y \big\}. 
$$
We then define a scaled version of $c$ by 
{
\begin{equation}  \label{e:def.cn.one.dim}
c_n(\Y, \Z; t) := c(r_n^{-1}\Y, r_n^{-1}\Z; t) =\one \big\{ \|y-z\| \ge r_nt, \ \text{for all } y \in \Y \text{ and } z\in \Z\setminus \Y \big\}. 
\end{equation}
}
and also, for a fixed $L\in (t,\infty)$, 
\begin{equation}  \label{e:def.sn}
s_n(\Y, \Z; t) := c_n(\Y, \Z; t)\, \one \big\{  \diam(\Y)\le r_n L \big\}. 
\end{equation}
Let $M_+\big( (\R^d)^k\big)$ be the space of Radon measures on $(\R^d)^k$.  {For a finite set $\Y\subset \R^d$ of $k$ points in general position, denote by $\ell(\Y)$ a \emph{center point} of $\Y$. For example, it may represent the left most point of $\Y$ in the lexicographic ordering. Another way of defining it is that one may set $\ell(\Y)$ to be 
a center of the unique $(k-2)$-dimensional sphere containing $\Y$. In either case, we  write $\overline \Y := \Y -\ell(\Y)$.}

The primary objective of this section is to describe the LDP for the point process 
\begin{equation}  \label{e:def.Lkn}
	\xi_{k,n} := \f1{\rho_{k,n}}\sum_{\Y \subset \Pn, \, |\Y|=k} s_n(\Y, \Pn; t) \, \delta_{r_n^{-1}\overline \Y}, \ \ n \ge 1. 
\end{equation}
The  process $\xi_{k,n}$ counts scaled $k$-tuples $r_n^{-1}\overline \Y$, that are locally concentrated in the sense of $\diam(\Y)\le r_nL$, under an additional restriction that  $r_n^{-1}\Y$ {be distant at least} $t$ from all the remaining points of $r_n^{-1}\Pn$. Notice that  $\diam(\Y) \le r_nL$ indicates $\| r_n^{-1}\overline \Y \|\le C$ for some constant $C\in(0,\infty)$. Thus, one can fix  a compact subset $E$ of $(\R^d)^k$ such that the process \eqref{e:def.Lkn} can be viewed as an element of $M_+(E)$. Because of this restriction, 
$M_+(E)$ is now equivalent to the space of \emph{finite} measures on $E$.   Write $\M_+(E)$ for the Borel $\sigma$-field generated by weak topology on $M_+(E)$. 

The proofs of the main results below are deferred to Section \ref{sec:LDP.pp.general.relative.entropy.and.rate}. 

\begin{theorem}  \label{t:LDP.pp.general}
The sequence $(\xi_{k,n})_{n\ge1}$ satisfies an LDP in the weak topology with rate $\rho_{k,n}$ and rate function 
{
\begin{equation}  \label{e:rate.func.pp.general}
\begin{split}
	\Lambda_k^*(\rho) &:=\sup_{f\in C_b(E)}\hspace{-.1cm} \Big\{ \int_Ef(\bx)  \rho (\dif \bx) \\
	&\qquad \qquad -\f1{k!} \int_{(\R^d)^{k-1}}\big( e^{f( \overline{(\bfz_d,\by)} )} -1 \big)\, \one \big\{   \diam (\bfz_d,\by) \le L \big\}\dif \by \Big\}, \ \ \rho\in M_+(E),
\end{split}
\end{equation}
}
where {$\bfz_d=(0,\dots,0)\in \R^d$}, $\by \in (\R^d)^{k-1}$ and $C_b(E)$ is the collection of continuous and bounded functions on $E$. 
\end{theorem}

The rate function  \eqref{e:rate.func.pp.general} can be associated  to the notion of \emph{relative entropy}. Writing $\lambda_m$ for Lebesgue measure on $(\R^d)^m$, we define 
\begin{equation}  \label{e:def.tau.k}
\tau_k(A) = \f1{k!}\, \lambda_{k-1}\big\{ \by\in (\R^d)^{k-1}: \diam(\bfz_d,\by)\le L, \, \overline{(\bfz_d,\by)} \in A  \big\}, \ \ A\subset E. 
\end{equation}
For a measure $\rho\in M_+(E)$, 
\begin{equation}  \label{e:def.rela.entropy}
H_k(\rho|\tau_k) := \begin{cases}
\int_E \log \f{\dif \rho}{\dif \tau_k}(\bx)\rho(\dif \bx) -\rho(E) +\tau_k(E) & \text{ if } \rho\ll \tau_k, \\
\infty & \text{ otherwise, }
\end{cases}
\end{equation}
denotes the relative entropy of $\rho$ with respect to $\tau_k$. In the special case when $\rho$ and $\tau_k$ are probability measures, \eqref{e:def.rela.entropy} reduces to the relative entropy defined for the space of probability measures; see, for example,  \cite[Equ.~(6.2.8)]{dembo:zeitouni:1998}. In our setup, however,  $\rho$ and $\tau_k$ are not necessarily probability measures, so we need to extend the definition of relative entropy as in \eqref{e:def.rela.entropy}. 

\begin{proposition}  \label{p:relative.entropy.and.rate}
Under the setup of Theorem \ref{t:LDP.pp.general}, we have that 
\begin{equation*}  \label{e:Hk.Lambdak}
\Lambda_k^*(\rho) = H_k(\rho|\tau_k), \ \ \ \rho\in M_+(E). 
\end{equation*}
\end{proposition}

{Finally, we also prove the analog of Theorem \ref{t:LDP.pp.general} when the Poisson point process is replaced by a binomial point process. To make this precise, we put }
$$
\xi_{k,n}^{\ms B} := \f1{\rho_{k,n}}\sum_{\Y \subset \Bn, \, |\Y|=k} s_n(\Y, \Bn; t) \, \delta_{r_n^{-1}\overline \Y}, \ \ n \ge 1, 
$$
{where $\Bn = \{X_1, \dots, X_n\}$ is a binomial point process consisting of $n$ i.i.d.~uniform random vectors in $[0, 1]^d$.}

\begin{corollary}  \label{c:LDP.pp.general}
The sequence $(\xi_{k,n}^{\ms B})_{n\ge1}$ satisfies an LDP in the weak topology with rate $\rho_{k,n}$ and rate function $\Lambda_k^*=H_k(\cdot | \tau_k)$.
\end{corollary}

%
%
\section{{Large deviation principles for geometric and topological functionals}}   \label{sec:application:geom.topo.funcl}

In this section, we provide the LDP for point processes {relating more directly to} geometric and topological {functionals. 
{Choose integers $k\ge 2$ and $m\ge1$, which remain fixed throughout the section}. Define 
 $H:=(h^{(1)},\dots,h^{(m)}): (\R^d)^k\to [0,\infty)^m$ to be a non-negative measurable function satisfying the following conditions: 
 \vspace{7pt}
 
\noindent \textbf{(H1)} $H$ is symmetric with respect to permutations of variables in $\R^d$. \\
\textbf{(H2)} $H$ is translation invariant:  
$$
H(x_1,\dots,x_k) = H(x_1+y, \dots, x_k+y),  \ \ \ x_i, y \in \R^d. 
$$
\textbf{(H3)} $H$ is locally determined: 
$$
H(x_1,\dots,x_k) = \bfz_m \ \ \text{whenever } \diam(x_1,\dots,x_k) > L, 
$$
where $L > 0 $ is a constant and $\bfz_m=(0,\dots,0)\in \R^m$. \\
\textbf{(H4)} For every  $\ba =(a_1,\dots,a_m)\in \R^m$, 
$$
\int_{(\R^d)^{k-1}} e^{\langle \ba, H(\bfz_d,\by)\rangle}  \Big( \one \big\{  H(\bfz_d,\by) \neq \bfz_m \big\} + \sum_{1 \le i \le  j \le m}h^{(i)}(\bfz_d,\by) h^{(j)}(\bfz_d,\by) \Big) \dif \by<\infty,
$$
where $\by\in (\R^d)^{k-1}$ and $\langle \cdot, \cdot \rangle$ denotes the Euclidean inner product.  \\
\textbf{(H5)} For every $\ba=(a_1,\dots,a_m)\in \R^m\setminus \{\bfz_m \}$, { it holds that $\int_{(\R^d)^{k-1}} \big| \sum_{i=1}^m a_i h^{(i)}(\bfz_d, \by) \big| \dif \by >0$.}
\vspace{7pt}

If $h^{(i)}$ are all {bounded functions},   we can drop condition \textbf{(H4)} because it can be implied by \textbf{(H3)}. Moreover,  note that \textbf{(H3)} remains true even when $L$ is increased. Hence, we may  assume that $L$ is larger than a fixed value $t > 0$. \\
Taking up a sequence $(r_n)_{n\ge1}$ in \eqref{e:sparse.radius}, we define a scaled version of $H$ by
\begin{equation}  \label{e:def.Hn}
H_n(x_1,\dots,x_k):= H(r_n^{-1}x_1, \dots, r_n^{-1}x_k), \ \ \ n=1,2,\dots. 
\end{equation}
For a $k$-point subset $\Y\subset \R^d$, a finite subset $\Z \supset \Y$ in $\R^d$ and $\bt=(t_1,\dots,t_m)\in [0,\infty)^m$, define 
\begin{equation}  \label{e:def.c}
c(\Y, \Z; \bt)  :=  \big( c(\Y, \Z; t_i) \big)_{i=1}^m, 
\end{equation}
and also, 
$$
G(\Y, \Z; \bt) := H(\Y)\odot c(\Y, \Z; \bt), 
$$
where $\odot$ means the Hadamard product; that is, for $A=(a_{ij})$ and $B=(b_{ij})$, we have $A\odot B=(a_{ij}b_{ij})$. We can define the scaled version of these functions by 
\begin{equation}  \label{e:def.cn.and.Gn}
c_n(\Y, \Z; \bt) := c(r_n^{-1}\Y, r_n^{-1}\Z; \bt), \ \  \text{ and } \ \ G_n(\Y, \Z; \bt) := H_n(\Y) \odot c_n(\Y, \Z; \bt). 
\end{equation}
For notational convenience, the $i$th entries of $H_n(\Y)$ and $G_n(\Y, \Z; \bt)$ are denoted respectively as $h_n^{(i)}(\Y)$ and  $g_n^{(i)}(\Y, \Z; t_i)$.  

Henceforth, our aim is to explore the large deviations of {the  process}
\begin{equation}  \label{e:def.Ukn}
U_{k,n} := \f1{\rho_{k,n}} \sum_{\Y \subset \Pn, \, |\Y|=k} \delta_{G_n(\Y, \Pn; \bt)} \in M_+(E'), 
\end{equation}
where $E'=[0,\infty)^m \setminus \{ \bfz_m\}$. 
Although the process \eqref{e:def.Ukn} aggregates the contributions of a certain score function over all $k$-tuples, it is not a pure $U$-statistic since this score function {depends also on} the remaining points of $\Pn$. 
Finally, as an analog of \eqref{e:def.tau.k}, we define 
\begin{equation*}  \label{e:def.tau.k.prime}
\tau_k'(A) := \f1{k!}\, \lambda_{k-1}\big\{ \by\in (\R^d)^{k-1}: H(\bfz_d,\by) \in A \big\}, \ \ \ A\subset E'. 
\end{equation*}
The proofs of all the theorems and propositions below are deferred to Sections \ref{sec:LDP.pp.Ustat} and \ref{sec:LDP.Ustat.and.rate.function}.

\begin{theorem}  \label{t:LDP.pp.Ustat}
The sequence $(U_{k,n})_{n\ge1}$ satisfies an LDP in the weak topology with rate $\rho_{k,n}$ and rate function 
\begin{equation}  \label{e:rate.func.pp.Ustat}
\begin{split}
	\bar \Lambda_k^*(\rho) &:=\sup_{f\in C_b(E')} \Big\{ \int_{E'} f(\bx)  \rho (\dif \bx) \\
	&\qquad \qquad -\f1{k!} \int_{(\R^d)^{k-1}} \big( e^{f( H(\bfz_d,\by) )} -1 \big)\, \one \big\{  H(\bfz_d,\by) \neq \bfz_m \big\} \dif \by \Big\}, \ \ \rho\in M_+(E'), 
\end{split}
\end{equation}
where $\by\in (\R^d)^{k-1}$. 
Furthermore, 
\begin{equation}  \label{e:rate.function.RE.pp.Ustat}
\bar \Lambda_k^*(\rho) = H_k' (\rho | \tau_k'), \ \ \rho\in M_+(E'), 
\end{equation}
where $H_k'$ is the relative entropy in $M_+(E')$ with respect to $\tau_k'$; that is, if $\rho\ll \tau_k'$, 
\begin{equation}  \label{e:def.Hk'}
H_k'(\rho|\tau_k') := \int_{E'}\log  \f{\dif \rho}{\dif \tau_k'}(\bx)\rho(\dif \bx) -\rho(E') +\tau_k'(E')
\end{equation}
and $H_k'(\rho|\tau_k') = \infty$ otherwise. 
\end{theorem}

{As in Corollary \ref{c:LDP.pp.general}, we can extend Theorem \ref{t:LDP.pp.Ustat} to the case of a binomial input. We define }
$$
U_{k,n}^{\ms B} := \f1{\rho_{k,n}} \sum_{\Y \subset \Bn, \, |\Y|=k} \delta_{G_n(\Y, \Bn; \bt)}. 
$$

\begin{corollary}  \label{c:LDP.pp.Ustat}
The sequence $(U_{k,n}^{\ms B})_{n\ge1}$ satisfies an LDP in the weak topology with rate $\rho_{k,n}$ and rate function $\bar \Lambda_k^*=H_k'(\cdot | \tau_k')$.
\end{corollary}

Subsequently, we also consider  statistics of the form
\begin{equation}  \label{e:def.Tkn}
T_{k,n}:= \sum_{\Y\subset \Pn, \, |\Y|=k} G_n(\Y, \Pn; \bt). 
\end{equation}
One can describe the corresponding LDP as follows. {We define $\mathcal D=\{ A\subset E': \bar A = \overline{\text{int}(A)} \}$.}

\begin{theorem}  \label{t:LDP.Ustat}
For every measurable {$A\in \mathcal D$}, we have as $n\to\infty$, 
\begin{equation}  \label{e:LDP.and.continuity}
\f1{\rho_{k,n}} \log \P \big( \rho_{k,n}^{-1} T_{k,n} \in A  \big) \to -\inf_{\bx\in A} I_k(\bx), 
\end{equation}
where $I_k$ is a rate function given as 
\begin{equation}  \label{e:def.Ik}
I_k(\bx) = \sup_{\ba\in \R^m} \Big\{ \langle \ba, \bx \rangle - \f1{k!}\, \int_{(\R^d)^{k-1}} \big( e^{\langle \ba, H(\bfz_d,\by)\rangle} - 1  \big)  \dif \by\Big\}, \ \ { \bx=(x_1,\dots,x_m)\in \R^m}, 
\end{equation}
with $\by\in (\R^d)^{k-1}$. 
\end{theorem}

{Again, we can extend the above result to the case of binomial point processes. We define }
$$
T_{k,n}^{\ms B}:= \sum_{\Y\subset \Bn, \, |\Y|=k} G_n(\Y, \Bn; \bt). 
$$

\begin{corollary}  \label{c:LDP.Ustat}
For every measurable {$A\in \mathcal D$}, we have as $n\to\infty$, 
	\begin{align*}
		\f1{\rho_{k,n}} \log \P \big( \rho_{k,n}^{-1} T_{k,n}^{\ms B} \in A  \big) \to -\inf_{\bx\in A} I_k(\bx). 
	\end{align*}
\end{corollary}

The rate function $I_k$ satisfies the following properties. 

\begin{proposition} \label{p:rate.function} 
$(i)$ $I_k$ is {continuously differentiable} on $E'$. \\
$(ii)$ $I_k$ is strictly convex on $E'$.  \\
$(iii)$ $I_k(\bx)=0$ if and only if $\bx=(\mu_1,\dots,\mu_m)$, where 
$$
\mu_i := \f1{k!}\, \int_{(\R^d)^{k-1}} h^{(i)}(\bfz_d,\by) \dif \by, \ \ \ i=1,\dots, m, 
$$
with $\by\in (\R^d)^{k-1}$. 
\end{proposition}

\begin{remark}   \label{rem:explicit.rate.function}
If $m=1$ and $h^{(1)}$ is an indicator function, {one can explicitly solve the variational problem in \eqref{e:def.Ik}, to get that}
 $I_k(x)=x \log (x/\mu_1)-x+\mu_1$. In this case, $I_k$ coincides with a rate function in the LDP for $\big( n^{-1}\sum_{i=1}^n X_i \big)_{n\ge1}$, where $(X_i)$ are i.i.d.~Poisson random variables with mean $\mu_1$. 
\end{remark}

\begin{example}[\v{C}ech complex component counts]  \label{ex:complex.component.counts}

We consider an application to {the \v{C}ech complex component counts}. Let $\check{C}( \X,  r )$ be the \v{C}ech complex on a point set {$\X=\{ x_1,\dots,x_m \} \subset\R^d$} with connectivity radius $r>0$. {Namely,}
\begin{itemize}
\item {The $0$-simplices of $\check{C}(\X, r)$ are the points in $\X$.}
\item {The $p$-simplex $\{x_{i_0},\dots,x_{i_p}\} \subset \X$ with $1\le i_0 < \dots < i_p\le m$, belongs to  $\check{C}(\X, r)$ if $\bigcap_{\ell=0}^p B(x_{i_\ell}, r/2)\neq \emptyset$, where $B(x,r)$ is the $d$-dimensional closed ball of radius $r$ centered at $x\in \R^d$.}
\end{itemize}
We then explore the LDP for 
\begin{equation}  \label{e:def.Skn}
S_{k,n} := \Big( \sum_{\Y\subset \Pn, \, |\Y|=k+1} \one \big\{ \check{C}(\Y, r_nt_i) \cong \Gamma_i \big\} \, c_n(\Y, \Pn; t_i) \Big)_{i=1}^m, 
\end{equation}
where $\Gamma_i$ is a \emph{connected} \v{C}ech complex on $k+1$ vertices, and $\cong$ means   isomorphism between simplicial complexes, and $c_n$ is defined in \eqref{e:def.cn.one.dim}. {Assume that $\Gamma_i\not\cong \Gamma_j$ for $i\neq j$.}
Then, the $i$th entry of $S_{k,n}$ represents the {number of components isomorphic to $\Gamma_i$} in the \v{C}ech complex $\check{C}(\Pn, r_nt_i)$. 
In this setting, it is easy to check that the function $H=(h^{(1)}, \dots, h^{(m)})$ with 
$$
h^{(i)} (x_1,\dots,x_{k+1}) = \one \Big\{  \check{C} \big( \{x_1,\dots,x_{k+1}\}, t_i \big) \cong \Gamma_i \Big\}, \ \ \ (x_1,\dots,x_{k+1}) \in (\bbr^d)^{k+1}, 
$$
satisfies conditions {\textbf{(H1)}--\textbf{(H5)}}. 

{Finally, we also define $S_{k,n}^\ms{B}$ to be the statistics analogous to \eqref{e:def.Skn}, {that are generated} by the binomial point process $\Bn$. 
The following corollary can be obtained as a direct application of  the above results.}
\end{example}

\begin{corollary}  
Assume that $\rho_{k + 1, n} \to \infty$ and $nr_n^d\to0$ as $n\to\infty$. Then, for any measurable {$A\in \mathcal D$}, as $n\to\infty$, 
\begin{align*}
&\f1{\rho_{k + 1, n}} \log \P \big( \rho_{k+1, n}^{-1} S_{k,n} \in A  \big) \to -\inf_{\bx\in A} I_{k+1}(\bx), \\ 
&{\f1{\rho_{k + 1, n}} \log \P \big( \rho_{k+1, n}^{-1} S_{k,n}^\ms{B} \in A  \big) \to -\inf_{\bx\in A} I_{k+1}(\bx).} 
\end{align*}
Here, $I_{k+1}$ is the rate function from \eqref{e:def.Ik}. The unique minimizer of $I_{k+1}$ equals $(\mu_1,\dots,\mu_m)$, where 
$$
\mu_i = \f1{(k+1)!}\, \int_{(\R^d)^k} \one \big\{  \check{C} \big( \{\bfz_d,\by\}, t_i \big) \cong \Gamma_i \big\}\dif \by, \ \ i=1,\dots,m, 
$$
with $\by \in (\R^d)^k$. 
\end{corollary}

{
\begin{example}[Component counts in a random geometric graph]
{We next deal with} the component counts in a random geometric graph $G(\Pn, r_nt_i)$ with radius $r_nt_i$, $i=1,\dots,m$. Define the function $H=(h^{(1)}, \dots, h^{(m)})$ by 
$$
h^{(i)} (x_1,\dots,x_{k}) = \one \Big\{  G \big( \{x_1,\dots,x_{k}\}, t_i \big) \cong \Gamma_i \Big\}, \ \ \ (x_1,\dots,x_{k}) \in (\bbr^d)^k, 
$$
where $\Gamma_i$ is a \emph{connected} graph on $k$ vertices and $\cong$ denotes a graph isomorphism. Assume that $\Gamma_i\not\cong \Gamma_j$ for $i\neq j$. Then, the collection of component counts defined by 
$$
S_{k,n}:= \Big( \sum_{\Y\subset \Pn, \, |\Y|=k} \one \big\{ G(\Y, r_nt_i) \cong \Gamma_i \big\} \, c_n(\Y, \Pn; t_i) \Big)_{i=1}^m, 
$$
satisfies the following LDP. 
\end{example}
\begin{corollary}  
	Assume that $\rho_{k,n} \to \infty$ and $nr_n^d\to0$ as $n\to\infty$. Then, for any measurable {$A\in \mathcal D$}, as $n\to\infty$, 
\begin{align*}
&\f1{\rho_{k,n}} \log \P \big( \rho_{k,n}^{-1} S_{k,n} \in A  \big) \to -\inf_{\bx\in A} I_{k}(\bx), \\
&{\f1{\rho_{k,n}} \log \P \big( \rho_{k,n}^{-1} S_{k,n}^\ms{B} \in A  \big) \to -\inf_{\bx\in A} I_{k}(\bx).}
\end{align*}
Once again, $I_{k}$ is the rate function with its unique minimizer $(\mu_1,\dots,\mu_m)$ given by 
$$
\mu_i = \f1{k!}\, \int_{(\R^d)^{k-1}} \one \big\{  G \big( \{\bfz_d,\by\}, t_i \big) \cong \Gamma_i \big\}\dif \by, \ \ i=1,\dots,m, 
$$
where $\by \in (\R^d)^{k-1}$. 
\end{corollary}
}

%
%
\section{{Applications in stochastic topology}}  \label{sec:applications}

In this section, we elucidate how to apply our main results to derive LDPs for two key quantities in stochastic topology, namely  \emph{persistent Betti numbers} (Section \ref{alpha_sec}) and \emph{Morse critical points} (Section \ref{morse_sec}). {In both cases, we assume $d = 2$ and  $k = 3$ in the notation of Section \ref{sec:model.main.results}.  In other words, we work with 1-dimensional topological  features generated by $3$ random points in the plane. 
Although the proposed methods can also be applied to $0$-dimensional features in the plane, we focus only on $1$-dimensional quantities, because they are much more non-trivial than $0$-dimensional ones. 
We conjecture that our findings can be generalized, at least partially, to higher-dimensional cases. 
The remark after Theorem \ref{t:LDP.persistent.Betti} below explains why such extensions are challenging.}

\subsection{Persistent Betti number of alpha complex}
\label{alpha_sec}
First, we deal with the persistent Betti numbers in the alpha complex. For readers not familiar with algebraic topology, {we briefly  discuss  conceptual ideas behind persistent Betti numbers. We suggest \cite{carlsson:2009,ghrist:2008} as a good introductory reading, while a more rigorous coverage of algebraic topology is in \cite{munkres:1996}.} 

{The Betti numbers $(\beta_k)_{k \ge 0}$ are fundamental invariants of topological spaces counting the number of $k$-dimensional cycles (henceforth we call it $k$-cycle)  as the boundary of a $(k+1)$-dimensional body. 
In the $3$-dimensional space, $\beta_1$ and $\beta_2$ can be viewed as the number of loops and cavities, respectively. Figure \ref{ball_fig} illustrates a sphere in $\R^3$, which encompasses one central cavity; therefore, $\beta_2 = 1$. The Betti number $\beta_1$ of this sphere is zero; even if
we wind a closed loop around the sphere, the loop ultimately vanishes as it moves upward (or
downward) along the sphere until the pole. 
Figure \ref{ball_fig} also illustrates a torus in $\R^3$, for which there are two distinct non-contractible loops; therefore, $\beta_1 = 2$. Moreover, the torus has a cavity, meaning that $\beta_2=1$.} 

\begin{figure}[!t]
\begin{center}
\includegraphics[width=10cm]{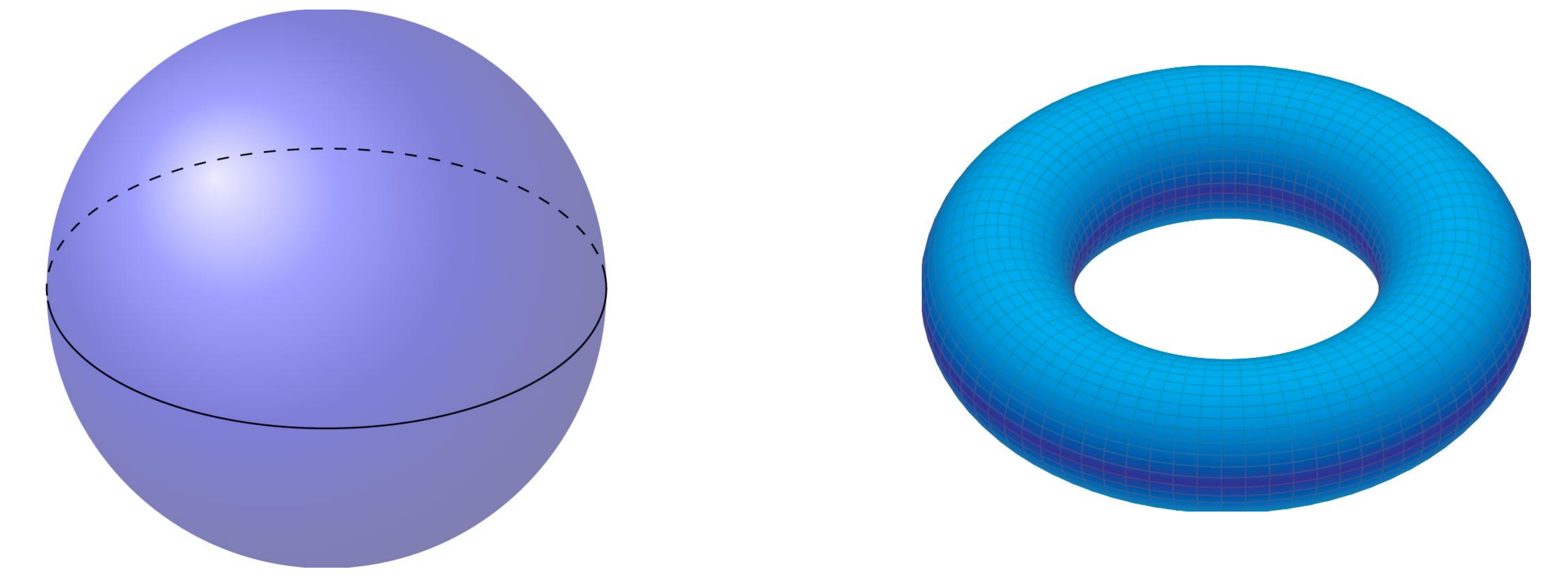}
\caption{{\footnotesize  Illustration of a sphere (left) and torus (right).}}
\label{ball_fig}
\end{center}
\end{figure}

We consider the case when the (persistent) Betti numbers are built on the alpha complex. The alpha complex is similar, in nature, to the \v{C}ech complex (see Example \ref{ex:complex.component.counts}), but it has a more natural geometric realization. Given a finite set $\X$ of points in $\R^d$ and $r\in[0,\infty]$, the alpha complex $\alpha (\X,r)$ is defined as a collection of subsets $\sigma \subset \X$ such that $\bigcap_{x\in \sigma} \big( B(x,r/2)\cap V_x \big)\neq \emptyset$, where $V_x$ is a Voronoi cell of $x$; that is, 
$V_x = \big\{ y\in \R^d: \| y-x \| \le \inf_{z \in \X} \| y-z\|\}$. Clearly, it holds that $\alpha(\X,r) \subset \check{C}(\X,r)$. Moreover, we have inclusions $\alpha(\X,r)\subset \alpha (\X,r')$ for all $r\le r'$, which indicates that $\alpha(\X,r)$ is a subcomplex of the Delaunay complex $\text{Del}(\X):=\alpha (\X,\infty)$. This property plays a crucial role in our analysis; see Remark \ref{rem:extension.higher.dim}. 

{Returning to the setup in Section \ref{sec:model.main.results}}, we consider the filtration induced by a collection of alpha complexes over a scaled Poisson point process $r_n^{-1} \Pn$, 
\begin{equation}  \label{e:def.alpha.complex}
\big( \alpha (r_n^{-1}\Pn, t), \, t\ge 0 \big) = \big( \alpha (\Pn, r_nt), \, t\ge 0 \big). 
\end{equation}
Topological changes in the filtration \eqref{e:def.alpha.complex} can be captured by the $k$th persistent Betti number $\beta_{k,n}(s,t)$ {for $0\le s \le t \le \infty$}.
Loosely speaking, $\beta_{k,n}(s, t)$ {represents the number of $k$-cycles}, that  appear in \eqref{e:def.alpha.complex} before time $r_ns$ and remain alive at time $r_nt$. More formally, $\beta_{k,n}(s,t)$ is defined as  
$$
\beta_{k,n}(s,t)  =\text{dim} \f{Z_k\big( \alpha(\Pn, r_ns) \big)}{Z_k\big( \alpha(\Pn, r_ns) \big) \cap B_k\big( \alpha(\Pn, r_nt) \big)}, 
$$
where $Z_k(\cdot)$ is the $k$th cycle group of an alpha complex and $B_k(\cdot)$ denotes its $k$th boundary group. Here, the homology coefficients are taken from an arbitrary field.  In the special case $s=t$, $\beta_{k,n}(s,t)$ reduces to the ordinary $k$th Betti number. 

{As mentioned before, in the following we restrict ourselves to a lower-dimensional case $d=2$ and $k=3$ (in the notation of Section \ref{sec:model.main.results}).}
For $(x_1, x_2,x_3)\in (\R^2)^3$ and $r>0$, define 
\begin{align}
h_r(x_1, x_2, &x_3) :=\one \Big\{ \beta_1\big( \alpha (\{x_1, x_2,x_3\}, r)  \big)=1 \Big\} \label{e:def.h.persistent.betti}      \\
&= \one \bigg\{  \Big\{  \bigcap_{j=1, \, j \neq j_0}^3 B(x_j, r/2) \neq \emptyset \ \text{for all } j_0\in \{ 1,2,3 \}\Big\} \cap \Big\{  \bigcap_{j=1}^3 B(x_j, r/2) =\emptyset\Big\}   \bigg\}. \notag 
\end{align}
For a  subset $\Y$ of $3$ points in $\R^2$, a finite set $\Z\supset \Y$ in $\R^2$, and $0\le s \le t <\infty$, 
\begin{align*}
g_{r_ns, r_nt}(\Y, \Z) &:= h_{r_ns}(\Y)\, h_{r_nt}(\Y)\, \one \big\{  \| y-z \| \ge r_nt,  \text{ for all } y\in \Y \text{ and } z\in \Z\setminus \Y  \big\}\\
&=h_{r_ns}(\Y)\, h_{r_nt}(\Y)\, \one \big\{  \alpha (\Y,r_nt) \text{ is a connected component of } \alpha(\Z, r_nt) \big\}. 
\end{align*}
Note that $g_{r_ns, r_nt}(\Y, \Z) =1$ if and only if a set $\Y$ in $\R^2$ with $|\Y|=3$ forms a single $1$-cycle before time $r_ns$, such that  this $1$-cycle remains alive at time  $r_nt$ and isolated from  all the remaining points in $\Z$ at that time. The proof of the next theorem is presented in Section \ref{sec:LDP.critical.point.and.LDP.persistent.Betti}. 
\begin{theorem}  \label{t:LDP.persistent.Betti}
Assume that $\rho_{3, n}\to\infty$ and $nr_n^2\to0$ as $n\to\infty$. Then, for every $0\le s_i \le t_i <\infty$, $i=1,\dots,m$, with $(s_i,t_i)\neq (s_j, t_j)$ for $i\neq j$, and {$A\in \mathcal D$}, 
\begin{equation}  \label{e:LDP.persistent.Betti.number}
\f1{\rho_{3, n}} \log \P \Big( \big(\rho_{3, n}^{-1}  \beta_{1,n}(s_i,t_i), \, i=1,\dots,m \big)\in A \Big) \to -\inf_{\bx\in A}I_3(\bx), \ \ \text{as } n\to\infty, 
\end{equation}
where $I_3$ is the rate function from \eqref{e:def.Ik}. The unique minimizer of $I_3$ equals $(\mu_1,\dots,\mu_m)$, where 
$$
\mu_i = \f1{6}\, \int_{(\R^2)^{2}} h_{s_i}(\bfz_2,\by)\, h_{t_i}(\bfz_2,\by) \dif \by, \ \ i=1,\dots,m, 
$$
with $\bfz_2=(0,0)\in \R^2$ and $\by\in (\R^2)^{2}$. \\
{Moreover, define $\beta_{1,n}^\ms{B}(s,t)$ to be the first-order persistent Betti number generated by the binomial point process. Then, $(\beta_{1,n}^\ms{B})_{n\ge1}$ satisfies the same LDP as \eqref{e:LDP.persistent.Betti.number}.}
\end{theorem}

{
\begin{remark}[Extensions to higher-dimensions]  \label{rem:extension.higher.dim}
In the proof of Theorem \ref{t:LDP.persistent.Betti}, we shall exploit the Morse inequality, which enables to bound the first-order persistent Betti number by the number of $1$-simplex counts (i.e., edge counts). A key property for dealing with the corresponding exponential moments is that the $1$-simplex count in the planar Delaunay triangulation grows at most linearly in the number of vertices (see Equ.~\eqref{e:d=2.k=3.Betti} for details). Unfortunately, this property breaks down in higher-dimensions (\cite{moment}); this is the reason why our discussion must be restricted to a special case $d=2$, $k=3$. 
Notice, however, that one needs this assumption only for showing \eqref{e:d=2.k=3.Betti}; the rest of our analyses holds true for general $d$ and $k$. 
For the LDP of higher-order persistent Betti numbers in a higher-dimensional space, we need to develop a new machinery that does not rely on the Morse inequality. More concretely, we need to detect the parameters $d$ and $k$, such that the (persistent) Betti number  grows at most linearly in the number of vertices. 
\end{remark}
}

%
%
\subsection{Morse critical points of min-type distance functions}
\label{morse_sec}

The objective of this example is to deduce the LDP for the number of Morse critical points of a certain min-type distance function. The behavior of Morse critical points of such distance functions have  been intensively investigated  in the context of central limit theorems and  the Poisson process approximation \cite{bobrowski:adler:2014, bobrowski:mukherjee:2015, bobrowski:schulte:yogeshwaran:2021,yogeshwaran:adler:2015}. 
{In addition to its intrinsic interests, this concept has served as a  practical quantifier of homological changes in random geometric complexes, especially in the field of Topological Data Analysis (\cite{bobrowski:adler:2014, bobrowski:kahle:2018}).}

Once again, we treat only the special case $d=2$ and $k=3$ (in the notation of Section \ref{sec:model.main.results}).
Given a homogeneous Poisson point process $\Pn$ on $[0, 1]^2$ {with intensity $n$}, we define min-type distance functions by 
$$
d_{\Pn}(x) := \min_{y\in \Pn} \| x-y \|, \ \ \ x\in \R^2. 
$$
Though $d_{\Pn}$ is not differentiable, one can still  define the notion of Morse critical points in the following sense. A point $c\in \R^2$ is said to be a Morse critical point of $d_{\Pn}$  with index $2$ if there exists a subset $\Y\subset \Pn$ of three points such that 
\vspace{5pt}

\noindent $(i)$ The points in $\Y$ are in general position. \\
$(ii)$ $d_{\Pn}(c)=\|c-y\|$ for all $y\in \Y$ and $d_{\Pn}(c)<\min_{z\in \Pn\setminus \Y}\|c-z\|$.\\
$(iii)$ The interior of the convex hull spanned by the points in $\Y$, denoted $\text{conv}^\circ (\Y)$, contains $c$. 
\vspace{3pt}

By the Nerve theorem (see, e.g., Theorem 10.7 in \cite{bjorner:1995}), for each $r>0$ the sublevel set $d_{\Pn}(-\infty,r]$ is homotopy equivalent to a \v{C}ech complex $\C(\Pn,2r)$. 
For this reason, the number of critical points of $d_{\Pn}$ with index $2$ whose critical values are less than $r_n$, behaves very similarly to the first-order Betti number of  $\C(\Pn,2r_n)$ (see \cite{bobrowski:adler:2014}). {A similar analysis was conducted for the case that distributions are supported on a closed manifold embedded in the ambient Euclidean space \cite{bobrowski:mukherjee:2015}. Moreover, \cite{yogeshwaran:adler:2015} studied the case when random points are sampled from a stationary point process.} 

Given {a point set $\Y\subset \R^2$ with $|\Y|=3$ in general position, let $\gamma(\Y)$ be the center of a unique $1$-dimensional sphere containing $\Y$, and define $\mathcal R(\Y) =d_{\Pn}\big( \gamma(\Y) \big)$. If $\gamma(\Y)$ gives a critical point of $d_{\Pn}$ (with index $2$), then $\mathcal R(\Y)$ represents its critical value.} Additionally,  $\mathcal U(\Y)$ denotes an open ball in $\R^2$ with radius $\mathcal R(\Y)$ centered at $\gamma (\Y)$. {Given $t_i\in [0,\infty)$, $i=1,\dots,m$, with $t_i\neq t_j$ for $i\neq j$,} we define 
\begin{equation}  \label{e:def.Nkn}
N_{n}:= \Big( \sum_{\Y\subset \Pn, \, |\Y|=3} \one \big\{ \gamma (\Y)\in \text{conv}^\circ (\Y), \ \mathcal R(\Y)\le r_nt_i, \ \mathcal U(\Y) \cap \Pn= \emptyset\big\} \Big)_{i=1}^m. 
\end{equation}
In particular, the $i$th entry of $N_{n}$ represents the number of Morse critical points of index $2$ with critical values less than or equal to $r_nt_i$.

The result below is essentially a  consequence of Theorem \ref{t:LDP.Ustat} and Proposition \ref{p:rate.function}. {Unlike Theorem \ref{t:LDP.persistent.Betti}, we do not provide the result for the version of binomial point processes; the required extension actually involves more complicated machinery. A formal proof is given in Section \ref{sec:LDP.critical.point.and.LDP.persistent.Betti}. 

\begin{theorem}  \label{t:LDP.critical.point}
Assume that $\rho_{3, n}\to\infty$ and $nr_n^2\to 0$ as $n\to\infty$. Then, for every measurable {$A\in \mathcal D$}, 
$$
	\f1{\rho_{3, n}}\, \log \P \big( \rho_{3, n}^{-1}N_{n}\in A \big) \to -\inf_{\bx\in A} I_3(\bx), \ \ \text{as } n\to\infty, 
$$
where $I_3$ is the rate function from \eqref{e:def.Ik}. The unique minimizer of $I_3$ equals $(\mu_1,\dots,\mu_m)$, where 
$$
\mu_i = \f1{6}\, \int_{(\R^2)^2} \one \big\{  \gamma (\bfz_2,\by) \in \text{conv}^\circ (\bfz_2,\by), \, \mathcal R(\bfz_2,\by) \le t_i \big\}\dif \by, \ \ i=1,\dots,m, 
$$
with $\by\in (\R^2)^2$. 
\end{theorem}

\section{Proofs}  \label{sec:proofs}

As preparation, we  partition the unit cube $[0,1]^d$ into $\rho_{k,n}$ sub-cubes $Q_1, \dots, Q_{\rho_{k,n}}$ of volume $\rho_{k,n}^{-1}$. To avoid notational complication we assume that $\rho_{k,n}$ takes only positive integers for all $n\in \bbn$. This assumption  applies to many of the sequences and functions throughout {this section}.   In particular, we set $Q_1=[0,\rho_{k,n}^{-1/d}]^d$. {Given a finite set $\Y$ of $k$ points in $\R^d$, we take  $\ell(\Y)$ to be the left most point of $\Y$ in the lexicographic ordering. As discussed in Section \ref{sec:model.main.results}, one may set $\ell(\Y)$ to be a center of the unique $(k-2)$-dimensional sphere containing $\Y$. In this case, however,  the description will be slightly more involved. Hence, for ease of description, we prefer to define $\ell(\Y)$ as the left most point.}

In the below, $C^*$ denotes a generic constant, which is independent of $n$  and may vary between and within the lines. 

\subsection{Proofs of Theorem \ref{t:LDP.pp.general}, Proposition \ref{p:relative.entropy.and.rate}, and Corollary \ref{c:LDP.pp.general}} \label{sec:LDP.pp.general.relative.entropy.and.rate}

{The proof of Theorem \ref{t:LDP.pp.general} can be completed via Propositions \ref{p:LDP.leading.term} and \ref{p:exp.negligible} below.} 
Proposition \ref{p:LDP.leading.term} is aimed to give the LDP for 
\begin{equation}  \label{e:def.eta.kn}
\eta_{k,n}=\f1{\rho_{k,n}} \sum_{\ell=1}^{\rho_{k,n}} \sum_{\Y\subset \Pnr_{Q_\ell}, \, |\Y|=k} s_n(\Y, \Pnr_{Q_\ell}; t)\, \delta_{r_n^{-1}\overline\Y} =: \f1{\rho_{k,n}} \sum_{\ell=1}^{\rho_{k,n}}  \eta_{k,n}^{(\ell)} \in M_+(E),
\end{equation}
where $\Pnr_{Q_\ell}$ denotes the restriction of $\Pn$  to the cube $Q_\ell$. 
{Setting up a ``blocked" point process as in \eqref{e:def.eta.kn} is a standard approach in the literature (\cite{seppalainen:yukich:2001, schreiber:yukich:2005}). The process \eqref{e:def.eta.kn} is of course different from the original process $(\xi_{k,n})_{n\ge1}$. In fact, if  geometric/topological objects exist spreading over multiple cubes in $(Q_\ell)_{\ell\ge1}$, then these objects are counted possibly by $(\xi_{k,n})_{n\ge1}$, whereas they are never counted by $(\eta_{k,n})_{n\ge1}$. {Despite such a difference}, one may justify in Proposition \ref{p:exp.negligible} that the difference between $(\xi_{k,n})_{n\ge1}$ and $(\eta_{k,n})_{n\ge1}$ is exponentially negligible in terms of the total variation distance. 
By virtue of this proposition, as well as Theorem 4.2.13 in \cite{dembo:zeitouni:1998}, our task can be reduced to prove the LDP for $(\eta_{k,n})_{n\ge1}$.}

{The main benefit of working with $(\eta_{k,n})_{n\ge1}$ is that it breaks down into a collection of i.i.d.~point processes $(\eta_{k,n}^{(\ell)})_{\ell\ge1}$.  
Note that $(\eta_{k,n}^{(\ell)})_{\ell\ge1}$ are defined in the space $M_p(E)$ of \emph{finite}  point measures on $E$. We here equip $M_p(E)$ with a Borel $\sigma$-field $\M_p(E):= \mathcal M_+(E) \cap M_p(E)$. Because of this decomposition, one can clarify the correspondence between \eqref{e:def.eta.kn} and the process  
\begin{equation}  \label{e:def.zeta.kn}
\zeta_{k,n} := \f1{\rho_{k,n}} \sum_{\ell=1}^{\rho_{k,n}} \zeta_k^{(\ell)}\in M_+(E), \ \ \ n=1,2,\dots,
\end{equation}
where $(\zeta_k^{(\ell)})_{\ell\ge1}$ are i.i.d.~Poisson random measures on $E$ with intensity measure $\tau_k$. 
Write $(\Omega', \mathcal F', \P')$ for the probability space on which \eqref{e:def.zeta.kn} is defined, and $\E'$ denotes the corresponding expectation.}

{
In Lemma \ref{l:LDP.zeta.kn} below, we first prove that \eqref{e:def.zeta.kn} satisfies the required LDP in Theorem \ref{t:LDP.pp.general}. After that, Lemmas \ref{l:total.variation.eta.kn} and \ref{l:exp.equiv.maximal.coupling} are aimed to establish exponential equivalence between \eqref{e:def.eta.kn} and \eqref{e:def.zeta.kn}. Now, we  formally state one of the main results of this section.}

\begin{proposition}  \label{p:LDP.leading.term}
The sequence $(\eta_{k,n})_{n\ge1}$ fulfills an LDP in the weak topology with rate $\rho_{k,n}$ and the rate function $\Lambda_k^*$ in  \eqref{e:rate.func.pp.general}. 
\end{proposition}

{Since the proof of Proposition \ref{p:LDP.leading.term} is rather long, we divide it into several lemmas. Combining these lemmas can conclude Proposition \ref{p:LDP.leading.term}.}
\begin{lemma}  \label{l:LDP.zeta.kn}
The sequence $(\zeta_{k,n})_{n\ge1}$ satisfies an LDP in the weak topology with rate $\rho_{k,n}$ and rate function $\Lambda_k^*$. 
\end{lemma}

\begin{proof}
It follows from Theorem 5.1 in \cite{resnick:2007} that for every $f\in C_b(E)$, 
\begin{align*}
\Lambda_k(f) 
&:= \log \E' \Big[ e^{\zeta_k^{(1)} (f) } \Big] = \f1{k!} \int_{(\R^d)^{k-1}} \big( e^{f(\overline{(\bfz_d,\by)})} -1 \big)\, \one \big\{  \diam(\bfz_d,\by) \le L \big\}\dif \by. 
\end{align*}
According to Cram\'er's theorem in Polish spaces in \cite[Theorem 6.1.3]{dembo:zeitouni:1998} (see also \cite[Corollary 6.2.3]{dembo:zeitouni:1998}), it turns out that $(\zeta_{k,n})_{n\ge1}$ satisfies a \emph{weak} LDP in  $M_+(E)$ with  rate function $\Lambda_k^*$. 

In order to extend this to a \emph{full} LDP, we need to demonstrate that $(\zeta_{k,n})_{n\ge1}$ is exponentially tight in the space $M_+(E)$.  The proof  is analogous to that in Lemma 6.2.6 of \cite{dembo:zeitouni:1998}. Since $\tau_k$ is tight in the space $E$, for every $\ell\ge 1$ there exists a compact subset $\Gamma_\ell\subset E$, such that $\tau_k (E \setminus \Gamma_\ell) \le \ell / (e^{2 \ell^2}-1)$. Define 
$$
K_\ell := \big\{ \nu \in M_+(E): \nu(E\setminus \Gamma_\ell) \le \ell^{-1}, \, \nu(E)\le \ell \big\}. 
$$
By Portmanteau's theorem for weak convergence, one can deduce that $K_\ell$ is weakly closed in $M_+(E)$. For $m\ge 1$, let $L_m:=\bigcap_{\ell=m}^\infty K_\ell$. Obviously, $L_m$ is weakly closed in $M_+(E)$. Furthermore,  Prohorov's theorem (see, e.g., Theorem A2.4.I in \cite{daley:verejones:2003}) ensures that  $L_m$ is relatively compact. Since  $L_m$ is now compact in $M_+(E)$, the exponential tightness of $(\zeta_{k,n})_{n\ge1}$ follows from 
\begin{equation}  \label{e:exp.tight}
\limsup_{n\to\infty} \f1{\rho_{k,n}} \log \P' \big(\zeta_{k,n}\in M_+(E)\setminus L_m\big) \le -\big( m-\tau_k(E)(e-1) \big), 
\end{equation}
for every $m\ge1$. 
By the union bound, 
$$
\P' \big(\zeta_{k,n}\in M_+(E)\setminus L_m\big) \le \sum_{\ell=m}^\infty \Big\{ \P'\big( \zeta_{k,n}(E\setminus \Gamma_\ell) > \ell^{-1} \big) + \P'\big( \zeta_{k,n}(E) >\ell \big) \Big\}. 
$$
By Markov's inequality and the fact that $\zeta_k^{(1)}(E\setminus \Gamma_\ell)$ is Poisson distributed with mean $\tau_k(E\setminus \Gamma_\ell)$, 
\begin{align}
&\sum_{\ell=m}^\infty  \P'\big( \zeta_{k,n}(E\setminus \Gamma_\ell) > \ell^{-1} \big)  \le \sum_{\ell=m}^\infty e^{-2\rho_{k,n} \ell} \Big( \E' \big[ e^{2\ell^2 \zeta_k^{(1)}(E\setminus \Gamma_\ell) } \big] \Big)^{\rho_{k,n}} \label{e:union.bdd1}\\
&= \sum_{\ell=m}^\infty e^{-2\rho_{k,n} \ell +\rho_{k,n} \tau_k(E\setminus \Gamma_\ell)(e^{2\ell^2}-1) } \notag \le \sum_{\ell=m}^\infty e^{-\rho_{k,n} \ell} \le 2 e^{-\rho_{k,n} m}.  \notag
\end{align}
By the similar calculation, 
\begin{equation}  \label{e:union.bdd2}
\sum_{\ell=m}^\infty \P'\big( \zeta_{k,n}(E) >\ell \big) \le 2e^{-( m-\tau_k(E)(e-1) )\rho_{k,n}}. 
\end{equation}
Putting \eqref{e:union.bdd1} and \eqref{e:union.bdd2} together gives \eqref{e:exp.tight}, as desired. 
\end{proof}
\vspace{5pt}

By virtue of Lemma \ref{l:LDP.zeta.kn}, our goal is now to show that  $(\eta_{k,n})_{n\ge1}$ in \eqref{e:def.eta.kn} satisfies the same LDP as $(\zeta_{k,n})_{n\ge1}$. As a first step, we claim that   for every $\ell\ge1$, the total variation distance between the laws of $\eta_{k,n}^{(\ell)}$ and $\zeta_k^{(\ell)}$ tends to $0$ as $n\to\infty$. Write $\mathcal L(\xi)$ for the probability law of a random element $\xi$.
\begin{lemma}  \label{l:total.variation.eta.kn}
For every $\ell\ge1$, 
\begin{equation}  \label{e:total.variation.eta.kn}
d_{\ms{TV}} \big( \mathcal L(\eta_{k,n}^{(\ell)}), \mathcal L(\zeta_k^{(\ell)}) \big) := \sup_{A\in \M_p(E)} \big| \P(\eta_{k,n}^{(\ell)}\in A) - \P'(\zeta_k^{(\ell)}\in A)  \big| \to 0, \ \ \text{as } n\to\infty. 
\end{equation}
\end{lemma}
\begin{remark}  \label{rem:critical.regime}
One of the  possible extensions of Lemma \ref{l:total.variation.eta.kn} seems to consider the case, for which  $(r_n)_{n\ge1}$ belongs to the critical regime, i.e., $nr_n^d\to c \in (0,\infty)$ as $n\to\infty$. In this case, however, the asymptotic behaviors of \eqref{e:def.eta.kn} and \eqref{e:def.weta.kn} below are essentially different. More concretely, in the critical regime \eqref{e:diff.total.variation} no longer holds, since the last term in \eqref{e:comp.involved} does not vanish as $n\to\infty$.  
\end{remark}
\begin{proof}
We begin with defining a sequence of random measures, 
\begin{equation}  \label{e:def.weta.kn}
\weta_{k,n}=\f1{\rho_{k,n}} \sum_{\ell=1}^{\rho_{k,n}} \sum_{\Y\subset \Pnr_{Q_\ell}, \, |\Y|=k} \one \big\{ \diam(\Y)\le r_nL  \big\}\, \delta_{r_n^{-1}\overline \Y} =: \f1{\rho_{k,n}} \sum_{\ell=1}^{\rho_{k,n}}  \weta_{k,n}^{(\ell)}, \ \ \ n\ge1. 
\end{equation}
Note that $(\weta_{k,n}^{(\ell)})_{\ell\ge1}$ are i.i.d.~point processes in the space $M_p(E)$. We first show that 
\begin{equation}  \label{e:total.variation.weta.kn}
d_{\ms{TV}} \big( \mathcal L(\weta_{k,n}^{(\ell)}), \mathcal L(\zeta_k^{(\ell)}) \big) \to0, \ \ \text{as } n\to\infty. 
\end{equation}
According to Theorem 3.1 in \cite{decreusefond:schulte:thaele:2016}, \eqref{e:total.variation.weta.kn} follows if one can verify that 
$$
\sup_{A: \text{Borel  in } E} \Big|  \E \big[\weta_{k,n}^{(\ell)}(A)\big] - \E'\big[ \zetal_k(A) \big]\Big| \to 0, \ \ \ n\to\infty, 
$$
and 
\begin{align}
\gamma_n&:= \max_{1\le p \le k-1} n^{2k-p} \int_{(Q_\ell)^{2k-p}} \one \big\{ \diam(x_1,\dots,x_k) \le r_nL \big\} \label{e:second.requirement.point.proc.conv}\\
&\qquad \qquad\qquad\qquad \times \one \big\{  \diam(x_1,\dots,x_p, x_{k+1}, \dots, x_{2k-p})  \le r_nL \big\} \dif \bx \to 0, \ \ n\to\infty. \notag
\end{align}
By the multivariate Mecke formula for Poisson point processes (see, e.g., Chapter 4 in \cite{last:penrose:2017}),
$$
\E \big[ \weta_{k,n}^{(\ell)}(A) \big] = \f{n^k}{k!}\, \int_{(Q_\ell)^k}\one \big\{  \diam(x_1,\dots,x_k) \le r_nL, \ r_n^{-1}\overline{(x_1,\dots,x_k)} \in A\big\} \dif \bx. 
$$
Performing the change of variables $x_i = x +r_n y_{i-1}$, $i=1,\dots,k$ (with $y_0\equiv 0$), 
{
\begin{align*}
\E \big[ \weta_{k,n}^{(\ell)}(A) \big] &= \f{\rho_{k,n}}{k!}\, \int_{Q_\ell}\int_{(\R^d)^{k-1}} \one \big\{ \diam(\bfz_d,\by) \le L, \, \overline{(\bfz_d,\by)} \in A\big\} \\
&\qquad \qquad \qquad \qquad \times \prod_{i=1}^{k-1} \one \{ x+r_n y_i\in Q_\ell \} \dif \by \dif x. 
\end{align*}
}
On the other hand, 
$$
\E'\big[ \zetal_k(A) \big] = \tau_k(A) =\f{\rho_{k,n}}{k!}\, \int_{Q_\ell}\int_{(\R^d)^{k-1}} \one \big\{ \diam(\bfz_d,\by)\le L, \ \overline{(\bfz_d,\by)} \in A\big\}  \dif \by \dif x. 
$$
Thus, as $n\to\infty$, 
\begin{align*}
\sup_A\Big|  \E \big[\weta_{k,n}^{(\ell)}(A)\big] - \E'\big[ \zetal_k(A) \big]\Big| &\le  \f{\rho_{k,n}}{k!}\, \int_{Q_\ell}\int_{(\R^d)^{k-1}} \one \big\{ \diam(\bfz_d,\by)  \le L\big\}  \\
&\qquad\qquad   \times  \Big( 1- \prod_{i=1}^{k-1}\one \{ x+r_n y_i \in Q_\ell \} \Big) \dif \by \dif x \to 0. 
\end{align*}
As for \eqref{e:second.requirement.point.proc.conv}, the change of variables $x_i=x+r_ny_{i-1}$, $i=1,\dots,2k-p$ (with $y_0\equiv 0$), yields that 
\begin{align*}
\gamma_n&:= \max_{1\le p \le k-1} n^{2k-p} r_n^{d(2k-p-1)} \int_{Q_\ell}\int_{(\R^d)^{2k-p-1}} \one \big\{ \diam(\bfz_d,y_1,\dots,y_{k-1}) \le L\big\}  \\
&\qquad  \times \one \big\{  \diam(\bfz_d,y_1,\dots,y_{p-1}, y_k,\dots, y_{2k-p-1}) \le L \big\} \prod_{i=1}^{2k-p-1} \one \{  x+r_n y_i\in Q_\ell\}\dif \by \dif x \\
&\le \max_{1\le p \le k-1} (nr_n^d)^{k-p} \int_{(\R^d)^{2k-p-1}} \prod_{i=1}^{2k-p-1} \one \big\{  \|y_i\| \le L\big\} \dif \by \to0, \ \ \ n\to\infty. 
\end{align*}
Thus, \eqref{e:total.variation.weta.kn} has been established, and it remains to verify that 
\begin{equation}  \label{e:diff.total.variation}
d_{\ms{TV}} \big( \mathcal L(\weta_{k,n}^{(\ell)}), \mathcal L(\eta_{k,n}^{(\ell)}) \big) \to 0, \ \ \ n\to\infty. 
\end{equation}
Noting that $\weta_{k,n}^{(\ell)}$ and $\eta_{k,n}^{(\ell)}$ are both  defined in the same probability space, we have 
$$
d_{\ms{TV}} \big( \mathcal L(\weta_{k,n}^{(\ell)}), \mathcal L(\eta_{k,n}^{(\ell)}) \big) =\f1{2}\sup_{\|f\|_\infty \le1} \Big|  \E\big[ f(\weta_{k,n}^{(\ell)}) \big] - \E\big[ f(\eta_{k,n}^{(\ell)}) \big] \Big|, 
$$
where the supremum is taken over all $f:M_p(E)\to\R$ {with $\|f\|_\infty:=\text{esssup}_{x\in M_p(E)}|f(x)| \le1$}. 
Whenever $f(\weta_{k,n}^{(\ell)}) - f(\eta_{k,n}^{(\ell)})$ is non-zero, there exists  a subset $\Y$ of $k$ points in $\Pnr_{Q_\ell}$,  such that $\diam(\Y)\le r_nL$ and $c_n(\Y, \Pnr_{Q_\ell}; t)=0$.  This implies that 
\begin{align*}
d_{\ms{TV}} \big( \mathcal L(\weta_{k,n}^{(\ell)}), \mathcal L(\eta_{k,n}^{(\ell)}) \big) &\le \f1{2}\sup_{\|f\|_\infty \le1} \E \bigg[  \, \big| \, f(\weta_{k,n}^{(\ell)})  -  f(\eta_{k,n}^{(\ell)}) \, \big| \\
&\qquad \times \one\bigg\{ \bigcup_{\Y\subset \Pnr_{Q_\ell}}  \Big\{  \diam(\Y)\le r_nL, \, c_n(\Y, \Pnr_{Q_\ell}; t)=0 \Big\} \bigg\}\bigg] \\
&\le \E \bigg[ \sum_{\Y\subset \Pnr_{Q_\ell}} \one \big\{ \diam(\Y)\le r_nL \big\}\, \big( 1-c_n(\Y, \Pnr_{Q_\ell}; t) \big) \bigg]. 
\end{align*}
Writing $\B\big( \{x_1,\dots,x_k\}; r \big) := \bigcup_{i=1}^k B(x_i, r)$ and appealing to the Mecke formula for Poisson point processes, the rightmost term above is equal to 
\begin{equation}  \label{e:comp.involved}
\begin{split}
&\f{n^k}{k!}\int_{(Q_\ell)^k} \one \big\{ \diam(x_1,\dots,x_k)\le r_nL\big\} \Big( 1-e^{-n\text{vol} \big(\B(\{ x_1,\dots,x_k \}; r_nt) \cap Q_\ell  \big)} \Big) \dif \bx  \\
&= \f{\rho_{k,n}}{k!}\int_{Q_\ell}\int_{(\R^d)^{k-1}} \one \big\{ \diam(\bfz_d,\by)\le L \big\}\prod_{i=1}^{k-1} \one \{ x+r_ny_i\in Q_\ell \}  \\
&\qquad \qquad \qquad \times \Big(1-e^{-n\text{vol} \big(\B(\{x, x+r_ny_1,\dots,x+r_ny_{k-1} \}; r_nt) \cap Q_\ell \big) }  \Big) \dif \by \dif x  \\
&\le \f1{k!}\int_{(\R^d)^{k-1}} \one \big\{ \diam(\bfz_d,\by)\le L \big\} \Big(1-e^{-nr_n^d \text{vol} \big(\B(\{\bfz_d,\by\}; t)  \big) }  \Big) \dif \by.
\end{split}
\end{equation}
Since $nr_n^d\to0$ as $n\to\infty$, the last integral tends to $0$ as $n\to\infty$. 
Now, \eqref{e:diff.total.variation} is obtained and the proof of \eqref{e:total.variation.eta.kn} is completed. 
\end{proof}
\vspace{5pt}

{By the standard argument on the maximal coupling} (see \cite[Lemma 4.32]{kallenberg}), for every $\ell\ge 1$, there exists a coupling $(\heta_{k,n}^{(\ell)}, \hzetal_k)$ on some probability space $(\hat{\Omega}_\ell, \hat{\mathcal F}_\ell, \hat \P_\ell)$, such that $\heta_{k,n}^{(\ell)} \stackrel{d}{=} \eta_{k,n}^{(\ell)}$, $\hzetal_k\stackrel{d}{=} \zetal_k$, and 
\begin{equation}  \label{e:maximal.coupling}
\hat \P_\ell \big(\heta_{k,n}^{(\ell)}\neq \hzetal_k \big) = d_{\ms{TV}} \big( \mathcal L(\eta_{k,n}^{(\ell)}), \mathcal L(\zeta_k^{(\ell)}) \big)  \to 0,  \ \ \ n\to\infty, 
\end{equation}
{where the last convergence follows from Lemma \ref{l:total.variation.eta.kn}.}
Define $\hat \Omega =\prod_{\ell=1}^\infty \hat \Omega_\ell$, $\hat{\mathcal F} = \bigotimes_{\ell=1}^\infty \hat{\mathcal F}_\ell$,  $\hat \P = \bigotimes_{\ell=1}^\infty \hat \P_\ell$, {and $\hat \E$ is the corresponding expectation}; then $(\heta_{k,n}^{(\ell)}, \hzetal_k)_{\ell \ge 1}$  
becomes a sequence of i.i.d.~random vectors under $\hat \P$. Define $\heta_{k,n}$ and $\hat \zeta_{k,n}$  analogously to \eqref{e:def.eta.kn} and \eqref{e:def.zeta.kn}. 

According to \cite[Theorem 4.2.13]{dembo:zeitouni:1998}, if one can show that  $(\heta_{k,n})_{n\ge1}$ and $(\hzeta_{k,n})_{n\ge1}$ are exponentially equivalent under a coupled probability $\hat \P$ (in terms of  the total variation distance), then it will be concluded that $(\eta_{k,n})_{n\ge1}$ and $(\zeta_{k,n})_{n\ge1}$ exhibit the same LDP. A more precise statement is given as follows.

\begin{lemma}   \label{l:exp.equiv.maximal.coupling}
For every $\delta>0$, 
\begin{equation}  \label{e:exp.neg1}
\lim_{n\to\infty} \f1{\rho_{k,n}} \log \hat \P \Big( d_{\ms{TV}}\big( \heta_{k,n}, \hzeta_{k,n} \big)\ge \delta \Big) =-\infty. 
\end{equation}
\end{lemma}

\begin{proof}
By Markov's inequality, we have, for any $a>0$, 
\begin{align*}
\hat \P \Big( d_{\ms{TV}}\big( \heta_{k,n}, \hzeta_{k,n} \big)\ge\delta\Big) &\le \hat \P \Big(\sum_{\ell=1}^{\rho_{k,n}} d_{\ms{TV}}  \big( \heta_{k,n}^{(\ell)}, \hzetal_k \big) \ge \delta \rho_{k,n}\Big) \\
	&\le  e^{-a\delta \rho_{k,n}} \Big(\hat \E \Big[ e^{a d_{\ms{TV}}  \big( \heta_{k,n}^{(1)}, \, \hat \zeta_k^{(1)} \big) }\Big] \Big)^{\rho_{k,n}}.
\end{align*}
Hence, \eqref{e:exp.neg1} follows if we can show that, for every $a>0$, 
\begin{equation}  \label{e:conv.E.hat}
\lim_{n\to\infty} \hat \E \Big[ e^{a d_{\ms{TV}}  \big( \heta_{k,n}^{(1)}, \, \hat \zeta_k^{(1)} \big) }\Big] =1. 
\end{equation}
Because of \eqref{e:maximal.coupling}, we get that $d_{\ms{TV}}\big( \heta_{k,n}^{(1)}, \, \hat \zeta_k^{(1)} \big) \stackrel{\hat \P}{\to}0$ as $n\to\infty$. It thus remains to demonstrate that 
\begin{equation}  \label{e:UI.E.hat}
\limsup_{n\to\infty}  \hat \E \Big[ e^{a d_{\ms{TV}}  \big( \heta_{k,n}^{(1)}, \, \hzeta_k^{(1)} \big) }\Big] <\infty, \ \ \text{ for all } a>0. 
\end{equation}
In fact, \eqref{e:UI.E.hat} ensures the uniform integrability for the convergence  \eqref{e:conv.E.hat}. 
By the Cauchy--Schwarz inequality, 
\begin{align*}
\hat \E \Big[ e^{a d_{\ms{TV}}  \big( \heta_{k,n}^{(1)}, \, \hat \zeta_k^{(1)} \big) }\Big]  &\le \Big\{  \E\Big[ e^{2a \eta_{k,n}^{(1)}(E)} \Big] \Big\}^{1/2}  \Big\{  \E'\Big[  e^{2a\zeta_k^{(1)} (E)}\Big] \Big\}^{1/2} \\
&=  \Big\{  \E \big[ e^{2a \sum_{\Y\subset \Pnr_{Q_1}} s_n(\Y, \Pnr_{Q_1}; t)} \big] \Big\}^{1/2} \Big\{ \exp \big\{  \tau_k(E) (e^{2a}-1) \big\} \Big\}^{1/2}. 
\end{align*}
Now we need to verify that for all $a>0$, 
\begin{equation}  \label{e:UI1}
\limsup_{n\to\infty} \E \Big[ e^{a \sum_{\Y\subset \Pnr_{Q_1}} s_n(\Y, \Pnr_{Q_1}; t)} \Big] <\infty. 
\end{equation}
For the proof of \eqref{e:UI1}, we consider the diluted family of  cubes 
\begin{equation}  \label{e:diluted.cubes}
G:=\big\{ 4Lr_n z + [0,tr_n/\sqrt{d}]^d \subset Q_1: z \in \bbz^d  \big\}
\end{equation}
{(recall that we have taken $L>t$).} Then, $Q_1$ can be  covered by at most $(4L\sqrt{d}/t)^d$ many translates of $G$. Denote these translates  as $G_1, \dots, G_{(4L\sqrt{d}/t)^d}$ (with $G_1\equiv G$). Let $b_n:=\rho_{k,n}^{-1}/(4Lr_n)^d$ denote  the number of cubes (of side length $tr_n/\sqrt{d}$) that are contained in $G$. 
As mentioned at the beginning of Section \ref{sec:proofs}, we assume, without loss of generality, that $(4L\sqrt{d}/t)^d$ and $b_n$ take only positive integers. 
Suppose $\Y$ is a set of $k$ points with  $\diam(\Y)\le r_nL$, then there is a unique $j\in \{ 1,\dots, (4L\sqrt{d}/t)^d \}$, so that {the left most point  $\ell(\Y)$} belongs to one of the cubes in $G_j$. Hence, 
{
\begin{equation}   \label{e:LMP.ineq}
\sum_{\Y\subset \Pnr_{Q_1}} s_n\big( \Y, \Pnr_{Q_1}; t\big) = \sum_{j=1}^{(4L\sqrt{d}/t)^d} \sum_{\Y\subset \Pnr_{Q_1}}  s_n\big( \Y, \Pnr_{Q_1}; t \big)  \one \Big\{ \ell(\Y)\in \bigcup_{J\in G_j}J \Big\}. 
\end{equation}  
}
Write $G_1=\{ J_1,\dots,J_{b_n} \}$ with $J_1=\big[ 0,tr_n/\sqrt{d} \big]^d$. By \eqref{e:LMP.ineq},  H\"older's inequality, and the spatial independence and homogeneity of $\Pn$,  we  need to demonstrate that 
\begin{equation}  \label{e:equiv.spatial.indep}
\limsup_{n\to\infty} \bigg( \E \Big[ e^{a\sum_{\Y\subset \Pnr_{Q_1}} s_n( \Y, \Pnr_{Q_1}; t) \one \{ \ell(\Y)\in J_1 \}  }  \Big]  \bigg)^{b_n}<\infty. 
\end{equation} 
{A key observation for the proof of \eqref{e:equiv.spatial.indep} is  that there is at most a single $k$-point subset $\Y\subset \Pnr_{Q_1}$ such that $c_n(\Y, \Pnr_{Q_1}; t)=1$ and $\ell(\Y)\in J_1$. In fact, if there are two distinct point sets $\Y, \Y' \subset \Pnr_{Q_1}$ with $\ell(\Y)\in J_1$ and  $\ell(\Y')\in J_1$, then it must be that $c_n(\Y, \Pnr_{Q_1}; t)=c_n(\Y', \Pnr_{Q_1}; t)=0$. It turns out from this observation that}
$$
\sum_{\Y\subset \Pnr_{Q_1}} s_n\big( \Y, \Pnr_{Q_1}; t \big) \one \big\{ \ell(\Y)\in J_1\big\} 
$$
is a $\{0,1\}$-valued random variable; hence, 
\begin{align*}
&\bigg( \E \Big[ e^{a\sum_{\Y\subset \Pnr_{Q_1}} s_n( \Y, \Pnr_{Q_1}; t ) \one \{ \ell(\Y)\in J_1 \}  }  \Big]  \bigg)^{b_n} \\
&= \bigg( 1+(e^a-1) \P\Big( \sum_{\Y\subset \Pnr_{Q_1}} s_n \big( \Y, \Pnr_{Q_1}; t \big) \one \big\{  \ell(\Y) \in J_1 \big\} =1\Big) \bigg)^{b_n} \\
&\le \bigg( 1+(e^a-1) \E \Big[ \sum_{\Y\subset \Pn}  \one \big\{\diam(\Y)\le r_nL, \,  \ell(\Y)\in J_1 \big\}  \Big] \bigg)^{b_n}. 
\end{align*}
Repeating the same calculations as before, it is not hard to see that
$$
\E \Big[ \sum_{\Y\subset \Pn} \one \big\{\diam(\Y)\le r_nL, \, \ell(\Y)\in J_1\big\} \Big] =C^*\big( \rho_{k,n}r_n^d+o(1) \big). 
$$
Now, we obtain 
$$
\bigg( 1+(e^a-1) \E \Big[ \sum_{\Y\subset \Pn}  \one \big\{\diam(\Y)\le r_nL, \, \ell(\Y)\in J_1 \big\}  \Big] \bigg)^{b_n} \to e^{C^*(e^a-1)/(4L)^d}<\infty, 
$$
as desired. 
\end{proof}
\vspace{5pt}

{Now that the proof of Proposition \ref{p:LDP.leading.term} has been completed, we next verify that the difference between $(\xi_{k,n})_{n\ge1}$ and $(\eta_{k,n})_{n\ge1}$ in terms of the total variation distance, is exponentially negligible.} 

\begin{proposition}  \label{p:exp.negligible}
For every $\delta>0$, 
\begin{equation}  \label{e:exp.negligible.L.kn.eta.kn}
\lim_{n\to\infty} \f1{\rho_{k,n}} \log \P \big( d_{\ms{TV}} ( \xi_{k,n}, \eta_{k,n} )   \ge \delta  \big) = -\infty. 
\end{equation}
\end{proposition}
\begin{proof}
For $\ell=1,\dots,\rho_{k,n}$, define a collection of points in $Q_\ell$ that are distance at most $r$ from the boundary of $Q_\ell$: 
$$
\bQ_\ell(r):= \big\{  x\in Q_\ell: \inf_{y\in \partial Q_\ell} \|x-y\| \le r\big\}, \  \ \ r>0. 
$$
For a subset $A\subset E$, we discuss two distinct cases for which a {$k$-point set $\Y\subset [0,1]^d$} with $\diam(\Y)\le r_nL$, makes different contributions to $\xi_{k,n}(A)$ and $\eta_{k,n}(A)$. The first case is that the  point set $\Y$ ``crosses a boundary" between two neighboring sub-cubes, i.e., there exist distinct $\ell_1$ and $\ell_2$ such that
$$
\Y\cap Q_{\ell_1} \neq \emptyset, \ \ \ \Y\cap Q_{\ell_2} \neq \emptyset. 
$$ 
Then, $\Y$ must be contained  in $\bigcup_{\ell=1}^{\rho_{k,n}} \bQ_\ell\big((L+t)r_n\big)$  {such that $\Y\cap \bigcup_{\ell=1}^{\rho_{k,n}} Q_\ell^\partial (tr_n)\neq \emptyset$}. Furthermore, $\Y$ may increase the value of $\xi_{k,n}(A)$, while  the value of $\eta_{k,n}(A)$ is unchanged. The second case  is that there are two neighboring sub-cubes $Q_{\ell_1}$ and $Q_{\ell_2}$, together with two point sets $\Y_i\subset Q_{\ell_i}$, $i=1,2$, of cardinality $k$  with $\diam(\Y_i)\le r_nL$, such that $\inf_{y_i \in \Y_i, \, i=1,2}\|y_1-y_2\|\le r_nt$. It then holds that $\Y_i\subset \bQ_{\ell_i}\big((L+t)r_n\big)$ {with $\Y_i \cap Q_\ell^\partial (tr_n)\neq \emptyset$ for $i=1,2$. Moreover,}  $\Y_i$ may  increase the value of   $\eta_{k,n}(A)$, but the  value of $\xi_{k,n}(A)$ is unchanged. Putting these observations together, we conclude that 
{
\begin{align}
d_{\ms{TV}}(\xi_{k,n}, \eta_{k,n}) &\le \f1{\rho_{k,n}}\sum_{\Y\subset \Pn} s_n \big(\Y, \Pnr_{\bigcup_{\ell=1}^{\rho_{k,n}} Q_{\ell}^\partial ((L+t)r_n)}; t \big) \one \Big\{ \Y \cap \bigcup_{\ell=1}^{\rho_{k,n}} Q_\ell^\partial (tr_n)\neq \emptyset  \Big\}   \label{e:disect.approx.two}  \\
&\qquad \qquad   + \f1{\rho_{k,n}}\sum_{\ell=1}^{\rho_{k,n}} \sum_{\Y\subset \Pn}  s_n \big(\Y, \Pnr_{\bQ_\ell ((L+t)r_n)}; t  \big) \one \big\{  \Y \cap Q_\ell^\partial (tr_n)\neq \emptyset \big\}. \notag
\end{align}
}
Now,  \eqref{e:exp.negligible.L.kn.eta.kn}  will follow if one can show that for every $\delta>0$, 
{
\begin{align}  
&\lim_{n\to\infty} \f1{\rho_{k,n}}\, \log \P \Big(  \sum_{\Y\subset \Pn} s_n \big(\Y, \Pnr_{\bigcup_{\ell=1}^{\rho_{k,n}} Q_{\ell}^\partial ((L+t)r_n)}; t \big) \label{e:exp.negligible1}\\
&\qquad \qquad \qquad\qquad \qquad \times \one \Big\{ \Y \cap \bigcup_{\ell=1}^{\rho_{k,n}} Q_\ell^\partial (tr_n)\neq \emptyset  \Big\}    \ge \delta \rho_{k,n} \Big) =-\infty, \notag
\end{align}
and 
\begin{equation}  \label{e:exp.negligible2}
\lim_{n\to\infty} \f1{\rho_{k,n}}\,\log  \P \Big(\sum_{\ell=1}^{\rho_{k,n}} \sum_{\Y\subset \Pn } s_n \big(\Y, \Pnr_{\bQ_\ell ((L+t)r_n)}; t  \big) \one \big\{\Y \cap Q_\ell^\partial (tr_n)\neq \emptyset  \big\} \ge \delta \rho_{k,n} \Big) =-\infty. 
\end{equation}
}
The proof techniques for \eqref{e:exp.negligible1} and \eqref{e:exp.negligible2} are similar, so we show only \eqref{e:exp.negligible1}. To begin, for each $1\le j \le d$,  denote a collection of ordered $j$-tuples by 
$$
\I_j = \big\{ \bell=(\ell_1,\dots,\ell_j): 1\le \ell_1 < \dots < \ell_j \le d \big\}. 
$$
{
For $\bell=(\ell_1,\dots,\ell_j)\in \I_j$, we define a collection of disjoint hyper-rectangles by 
\begin{align}
	J \hspace{-.1cm}&:=\hspace{-.1cm} \bigg\{  \Big(  \rho_{k,n}^{-1/d}z  +\big[0,\rho_{k,n}^{-1/d}\big]^{\ell_1-1}\times \big[-(L+t)r_n,(L+t)r_n\big]\times \big[0,\rho_{k,n}^{-1/d}\big]^{\ell_2-\ell_1-1}  \label{e:translation.hyper.rectangle} \\
&\qquad \times \big[-(L+t)r_n,(L+t)r_n\big] \times \big[0,\rho_{k,n}^{-1/d}\big]^{\ell_3-\ell_2-1} \times \dots \times  \big[-(L+t)r_n,(L+t)r_n\big] \notag  \\
&\qquad \times \big[0,\rho_{k,n}^{-1/d}\big]^{\ell_j-\ell_{j-1}-1}\hspace{-.2cm} \times \hspace{-.1cm}\big[-(L+t)r_n,(L+t)r_n\big] \hspace{-.1cm}\times\hspace{-.1cm} \big[0,\rho_{k,n}^{-1/d}\big]^{d-\ell_j} \Big) \cap [0,1]^d: z\in \bbz_+^d \bigg\}. \notag 
\end{align}
By construction, all the rectangles in $J$ {are contained in $\bigcup_{\ell=1}^{\rho_{k,n}} Q_{\ell}^\partial ((L+t)r_n)$}. Moreover, since the number of rectangles in $J$ is $\rho_{k,n}$, one can enumerate $J$ in a way that $J=(I_{p,n}^\bell, \, p=1,\dots,\rho_{k,n})$. 
}
Then,  one can offer the following bound:
{
\begin{align}
&\sum_{\Y\subset \Pn} s_n \big(\Y, \Pnr_{\bigcup_{\ell=1}^{\rho_{k,n}} Q_{\ell}^\partial ((L+t)r_n)}; t \big) \one \Big\{ \Y \cap \bigcup_{\ell=1}^{\rho_{k,n}} Q_\ell^\partial (tr_n)\neq \emptyset  \Big\} \label{e:disecting.bound} \\   &\le \sum_{j=1}^d \sum_{\bell \in \I_j} \sum_{p=1}^{\rho_{k,n}}\sum_{\Y\subset \Pnr_{I_{p,n}^\bell}} \hspace{-10pt} s_n \big( \Y, \Pnr_{I_{p,n}^\bell}; t \big). \notag 
\end{align}
}
Owing to this bound, it  remains to show that for every $j\in \{ 1,\dots,d \}$, $\bell\in \I_j$, and $\delta>0$, 
\begin{equation}  \label{e:upper.negligible}
\lim_{n\to\infty} \f1{\rho_{k,n}}\, \log \P\Big(  \sum_{p=1}^{\rho_{k,n}} \sum_{\Y\subset \Pnr_{I_{p,n}^\bell}} \hspace{-10pt}s_n\big(\Y, \Pnr_{I_{p,n}^\bell}; t \big) \ge \delta \rho_{k,n} \Big) = -\infty. 
\end{equation}
Since $( I_{p,n}^\bell, \, p=1,\dots,\rho_{k,n} )$  are disjoint, the spatial independence of $\Pn$ ensures that 
$$
\bigg(  \sum_{\Y\subset \Pnr_{I_{p,n}^\bell}} \hspace{-10pt}s_n\big(\Y, \Pnr_{I_{p,n}^\bell}; t \big), \, p=1,\dots,\rho_{k,n} \bigg)
$$
are i.i.d.~random variables. 
Hence, we have for every $a>0$, 
\begin{align*}
&\P\Big(  \sum_{p=1}^{\rho_{k,n}} \sum_{\Y\subset \Pnr_{I_{p,n}^\bell }} \hspace{-10pt}s_n\big(\Y, \Pnr_{I_{p,n}^\bell}; t \big) \ge \delta \rho_{k,n} \Big) \\
	&\quad\le e^{-a\delta \rho_{k,n}} \bigg(  \E\Big[  \exp\Big\{a\hspace{-1pt}\sum_{\Y\subset \Pnr_{I_{1,n}^\bell}} \hspace{-10pt}s_n\big(\Y, \Pnr_{I_{1,n}^\bell}; t \big)\Big\} \Big] \bigg)^{\rho_{k,n}}. 
\end{align*}
For the proof of \eqref{e:upper.negligible}, it is now sufficient to show that for every $a>0$, 
\begin{equation}  \label{e:UI.plus.conv}
\E \bigg[ \exp\Big\{  a\hspace{-5pt}\sum_{\Y\subset \Pnr_{I_{1,n}^\bell}} \hspace{-10pt}s_n\big(\Y, \Pnr_{I_{1,n}^\bell}; t \big)   \Big\}  \bigg]  \to 1, \ \ \text{as } n\to\infty. 
\end{equation}
By the same argument as that for \eqref{e:UI1}, one can see that for every $a>0$, 
$$
\limsup_{n\to\infty} \E \bigg[ \exp\Big\{  a \hspace{-5pt}\sum_{\Y\subset \Pnr_{I_{1,n}^\bell}} \hspace{-10pt}s_n\big(\Y, \Pnr_{I_{1,n}^\bell}; t\big)   \Big\}  \bigg]<\infty, 
$$
which  implies the required uniform integrability. Now, \eqref{e:UI.plus.conv} will follow if we can verify that 
\begin{equation}  \label{e:conv.in.prob.comp}
\sum_{\Y\subset \Pnr_{I_{1,n}^\bell }} \hspace{-10pt}s_n\big(\Y, \Pnr_{I_{1,n}^\bell}; t \big)   \stackrel{\P}{\to} 0, \ \ \text{as  } n\to\infty. 
\end{equation}
To show this, we have as $n\to\infty$, 
\begin{align*}
&\E\Big[ \sum_{\Y\subset \Pnr_{I_{1,n}^\bell}} \hspace{-10pt}s_n\big(\Y, \Pnr_{I_{1,n}^\bell}; t \big)  \Big] \le \E\Big[  \sum_{\Y\subset \Pnr_{I_{1,n}^\bell }} \hspace{-10pt}\one \big\{ \diam(\Y)\le r_nL \big\} \Big] \label{e:exp.evaluation.comp}\\
&\qquad \qquad = \f{n^k}{k!}\, \int_{( I_{1,n}^\bell )^k}\one \big\{ \diam(x_1,\dots,x_k)\le r_nL \big\} \dif \bx \notag  \\
&\qquad \qquad = \f{\rho_{k,n}}{k!}\, \int_{I_{1,n}^\bell}\int_{(\R^d)^{k-1}} \one\big\{ \diam(\bfz_d,\by)\le L\big\} \prod_{i=1}^{k-1}\one \big\{ x+r_n y_i \in I_{1,n}^\bell \big\}\dif \by \dif x \notag \\
&\qquad \qquad \le \f{\rho_{k,n}}{k!}\, \text{vol}( I_{1,n}^\bell) \int_{(\R^d)^{k-1}} \one \big\{ \diam(\bfz_d,\by)\le L \big\} \dif \by  =\mathcal O \big( (nr_n^d)^{k/d} \big)\to 0. \notag
\end{align*}
Hence \eqref{e:conv.in.prob.comp} has been established, as desired. 
\end{proof}

{Before concluding this subsection, we present the proof of Proposition \ref{p:relative.entropy.and.rate}.}

\begin{proof}[Proof of Proposition \ref{p:relative.entropy.and.rate}]
Our proof is closely related to Theorem 5.4 in \cite{rassoulagha:seppalainen:2015}. 
Our goal here is to show that,  for every $f\in C_b(E)$, 
\begin{equation}  \label{e:duality.prop}
\int_E \big( e^{f(\bx)}-1 \big)\tau_k(\dif \bx) = \sup_{\rho\in M_+(E)} \big\{ \langle \rho, f\rangle -H_k(\rho | \tau_k) \big\}, 
\end{equation}
where $\langle \rho, f \rangle :=\int_Ef(\bx)\rho(\dif \bx)$. 
Given $f\in C_b(E)$, we define $\dif \rho=e^f \dif \tau_k$. Then, $\rho\in M_+(E)$ with $\rho\ll \tau_k$. By \eqref{e:def.rela.entropy}, it is elementary to calculate that   $\langle \rho, f\rangle -H_k(\rho|\tau_k) = \int_E \big( e^{f(\bx)}-1 \big)\tau_k(\dif \bx)$, which has shown that the left hand side in \eqref{e:duality.prop} is bounded by the right hand side. 

Next, let us take  $\nu\in M_+(E)$ with $\nu\ll \tau_k$, say with density $\varphi$. Then, by Jensen's inequality,
$$
\int_E \log \bigg( \frac{e^{f(\bx)}}{\varphi(\bx)}\bigg)\, \f{\nu(d \bx)}{\nu(E)} \le\log \bigg( \int_E \frac{e^{f(\bx)}}{\varphi(\bx)} \f{\nu(d \bx)}{\nu(E)} \bigg)=  \log \bigg(\int_Ee^{f(\bx)} \f{\tau_k(d \bx)}{\nu(E)}\bigg).
$$
Hence, by  \eqref{e:def.rela.entropy} {and the elementary inequality: $1+\log x \le x$ for $x>0$}, 
$$
\langle \nu, f\rangle -H_k(\nu|\tau_k)  { \le  \nu(E) \bigg\{ 1+\log \bigg( \int_E  e^{f(\bx)} \f{\tau_k(d x)}{\nu(E)}\bigg)  \bigg\} -\tau_k(E)} \le \int_E \big( e^{f(\bx)}-1 \big)\tau_k(\dif \bx). 
$$
Now \eqref{e:duality.prop} is obtained, and the rest of the argument after \eqref{e:duality.prop} is essentially the same as Theorem 5.4 in \cite{rassoulagha:seppalainen:2015}, so we will omit it. 
\end{proof}

{Finally, we prove Corollary \ref{c:LDP.pp.general}.  We can deduce its assertion from Theorem \ref{t:LDP.pp.general} by showing that for every $\varepsilon_0 > 0$,}  
$$\lim_{n \to \infty}\rho_{k, n}^{-1}\log\P\big(d_{\ms{TV}}(\xi_{k, n}, \xi_{k,n}^{\ms B}) \ge \varepsilon_0\big) = -\infty.$$

\begin{proof}[Proof of Corollary \ref{c:LDP.pp.general}]
{Similarly as in the proof of Lemma \ref{l:exp.equiv.maximal.coupling}, we consider a family of diluted cubes. In the current setting, we take}
\begin{equation}  \label{e:diluted.cubes2}
	G := \big\{ 8Lr_n z + [0, 4Lr_n]^d \subset [0, 1]^d: z \in \bbz^d  \big\}.
\end{equation}
Then, $[0,1]^d$ can be covered by $2^d$ translates of $G$. 
As before, we write $G=\{ J_1,\dots,J_{b_n'} \}$ with $J_1 = \big[0, 4Lr_n\big]^d$,  where $b_n' := (8Lr_n)^{-d}$ denotes the number of cubes that are contained in $G$. {Since there are at most finitely many translates of $G$, it suffices to prove that as $n\to\infty$, 
\begin{align*}
&\rho_{k,n}^{-1} \log \P \bigg(  \sup_{A\subset E} \Big| \sum_{\Y\subset \Pn, \, |\Y|=k} s_n(\Y, \Pn; t) \one \Big\{\ell(\Y)\in \bigcup_{i=1}^{b_n'} J_i\Big\} \, \delta_{r_n^{-1} \overline \Y}(A) \\
&\qquad \qquad \qquad -\sum_{\Y\subset \Bn, \, |\Y|=k} s_n(\Y, \Bn; t) \one \Big\{\ell(\Y)\in \bigcup_{i=1}^{b_n'}J_i \Big\} \, \delta_{r_n^{-1} \overline \Y}(A) \Big| \ge \vep_0\rho_{k,n} \bigg) \to -\infty 
\end{align*}  
(recall that $\ell(\Y)$ is the left most point of $\Y$ in the lexicographic ordering). 
}

We say that $J_i$ is an \emph{$n$-bad cube} if one of the following events happens. 
\begin{itemize}
\item There exists a $k$-element subset $\Y \subset \Pn$ with $\ell(\Y)\in J_i$, $s_n(\Y, \Pn; t)=1$, such that $\Y \not\subset \Bn$ or $s_n\big(\Y, \Bn; t\big) = 0$ holds. 
\vspace{3pt}
\item There exists a $k$-element subset $\Y \subset \Bn$ with $\ell(\Y)\in J_i$, $s_n(\Y, \Bn; t)=1$, such that $\Y \not\subset \Pn$ or $s_n\big(\Y, \Pn; t\big) = 0$ holds. 
\end{itemize} 
{In this setting, the key observation is that there exists a constant $M>0$, such that 
\begin{align*}
&\sup_{A\subset E} \Big| \sum_{\Y\subset \Pn, \, |\Y|=k} s_n(\Y, \Pn; t) \one \Big\{\ell(\Y)\in \bigcup_{i=1}^{b_n'} J_i\Big\} \, \delta_{r_n^{-1} \overline \Y}(A) \\
&\qquad \qquad \qquad -\sum_{\Y\subset \Bn, \, |\Y|=k} s_n(\Y, \Bn; t) \one \Big\{\ell(\Y)\in \bigcup_{i=1}^{b_n'}J_i \Big\} \, \delta_{r_n^{-1} \overline \Y}(A) \Big| \\
&\le M \sum_{i=1}^{b_n'}\one \{ J_i \text{ is $n$-bad} \}. 
\end{align*}
}
Thus, it  suffices to show that for every $\varepsilon_0 > 0$, 
$$
\lim_{n \to \infty}\frac{1}{\rho_{k, n}}\log\P\Big(\sum_{i=1}^{b_n'}\one \{ J_i \text{ is $n$-bad} \} \ge \varepsilon_0 \rho_{k, n} \Big) = -\infty.
$$
For $0<\varepsilon \le1$, let $\Pn^{(\varepsilon)}$ be a homogeneous Poisson point process on $[0,1]^d$ with intensity $n\varepsilon$, independent of $\Pn$. Then, $\Pn^{(\varepsilon, \ms{a})}:= \Pn \cup \Pn^{(\varepsilon)}$ represents an \emph{augmented}  version of $\Pn$ with intensity $n(1+\varepsilon)$.  Moreover, $\Pn^{(\varepsilon, \ms{t})}$ denotes a \emph{thinned} version of $\Pn$, obtained by removing points with probability $\varepsilon$. Notice that $\Pn^{(\varepsilon, \ms{a})}\stackrel{d}{=}\mathcal P_{n(1+\varepsilon)}$ and $\Pn^{(\varepsilon, \ms{t})}\stackrel{d}{=}\mathcal P_{n(1-\varepsilon)}$. 
In this setting,  we introduce the event 
$$F_{n, \varepsilon} := \big\{\Pn^{(\varepsilon, \ms{t})} \subset \Bn \subset \Pn^{(\varepsilon, \ms{a})}\big\},$$
	and note by the Poisson concentration bound from \cite[Lemma 1.2]{penrose:2003},
\begin{align}
	\label{fne_eq}
	\limsup_{n \to \infty}\rho_{k, n}^{-1}\log\P\big(F_{n, \varepsilon}^c\big) \le -C^*\lim_{n\to\infty} \rho_{k, n}^{-1} n = -\infty.
\end{align}
The key advantage of the event $F_{n, \varepsilon}$ is that it allows to simplify the property of being $n$-bad. Indeed, if $J_i$ is $n$-bad and $F_{n, \varepsilon}$ holds,  then $J_i$ becomes \emph{$(n, \varepsilon)$-special}, in the sense that 
$$
	T \cap (\Pn^{(\varepsilon, \ms{a})}  \setminus \Pn^{(\varepsilon, \ms{t})} )   \ne \emptyset, \ \ \text{ and } \ \ \Pn^{(\varepsilon, \ms{a})} (T) \ge k, 
$$
where 
$$
T := \big\{ x\in \R^d: \inf_{y\in J_i} \|x-y \| \le 2r_nL \big\}, \ r >0. 
$$
Since we work with diluted cubes in \eqref{e:diluted.cubes2}, one can see that $\big( \one \{ J_i \text{ is } (n,\varepsilon)\text{-special} \} \big)_{i=1}^{b_n'}$ are i.i.d.~Bernoulli random variables. Thus, the number of $(n, \varepsilon)$-special cubes is a binomial random variable with $b_n'$ trials and success probability 
$p_{n, \varepsilon} := \P\big(J_1 \text{ is $(n, \varepsilon)$-special}\big).$
Then, one can bound $p_{n,\varepsilon}$ as follows: 
\begin{align*}
p_{n,\varepsilon} &\le \P \big( T\cap (\Pn^{(\varepsilon, \ms{a})} \setminus \Pn^{(\varepsilon, \ms{t})}) \neq \emptyset \, \big| \, \Pn^{(\varepsilon, \ms{a})} (T)=k\big) \P\big( \Pn^{(\varepsilon, \ms{a})} (T)=k \big) \\
&\qquad \qquad \qquad \qquad \qquad \qquad \qquad\qquad + \P\big( \Pn^{(\varepsilon, \ms{a})} (T)\ge k+1 \big) \\
&\le 2k\varepsilon  \P\big( \mathcal P_{n(1+\varepsilon)} (T)=k \big)  + \P\big( \mathcal P_{n(1+\varepsilon)} (T)\ge k+1 \big) \\
&\le C^* \big( \varepsilon (nr_n^d)^k + (nr_n^d)^{k+1} \big). 
\end{align*}
{In conclusion, if one takes sufficiently small $\vep>0$, then 
$$
b_n' p_{n,\vep} \le C^* b_n' \big( \varepsilon (nr_n^d)^k + (nr_n^d)^{k+1} \big) \le \vep_0 \rho_{k,n}
$$
for large $n$ enough. 
Therefore,  the binomial concentration inequality \cite[Lemma 1.1]{penrose:2003} gives that
\begin{align*}
\limsup_{n\to\infty} \frac{1}{\rho_{k,n}} \log \P \big( \text{Bin} (b_n', p_{n,\vep}) \ge \varepsilon_0 \rho_{k,n} \big) &\le -\lim_{n\to\infty} \frac{\varepsilon_0 }{2}\, \log \bigg\{ \frac{\varepsilon_0 \rho_{k,n}}{C^* b_n' ( \varepsilon (nr_n^d)^k + (nr_n^d)^{k+1})}  \bigg\} \\
&= -\frac{\varepsilon_0}{2} \log \Big\{ \frac{(8L)^d \varepsilon_0}{C^*\varepsilon} \Big\}. 
\end{align*}
}
The last term tends to $-\infty$ as $\varepsilon\to0$. 
Hence, combining this result with \eqref{fne_eq} concludes the proof of Corollary \ref{c:LDP.pp.general}.
\end{proof}

\subsection{Proofs of Theorem \ref{t:LDP.pp.Ustat} and Corollary \ref{c:LDP.pp.Ustat}}  \label{sec:LDP.pp.Ustat}
{First  we point out that the proof of  Corollary \ref{c:LDP.pp.Ustat} is almost identical to that of Corollary \ref{c:LDP.pp.general}, so we skip it here. 
For the proof of Theorem \ref{t:LDP.pp.Ustat}, we define}
$$
V_{k,n} := \frac{1}{\rho_{k,n}} \sum_{\Y\subset \Pn, \, |\Y|=k} s_n(\Y, \Pn; t)\, \delta_{H_n(\Y)} \in M_+(E'). 
$$
The  key step of our proof  is to show that the assertions  of Theorem \ref{t:LDP.pp.Ustat} still hold even when $(U_{k, n})_{n\ge1}$ is replaced by $(V_{k, n})_{n\ge1}$. To clarify our presentation, we  state this step as a separate proposition. 
\begin{proposition}  \label{p:LDP.Vkn}
The sequence $(V_{k,n})_{n\ge1}$ satisfies an LDP in the weak topology with rate $\rho_{k,n}$ and the rate function $\bar \Lambda_k^*$  defined in  \eqref{e:rate.func.pp.Ustat}.
\end{proposition}
\begin{proof}
The process $(V_{k,n})_{n\ge1}$ has structure very similar to that of $(\xi_{k,n})_{n\ge1}$ in \eqref{e:def.Lkn}; thus, the proof techniques for Proposition \ref{p:LDP.Vkn} are  parallel to those for Theorem \ref{t:LDP.pp.general}. More precisely, as an analog of \eqref{e:def.eta.kn}, we define 
$$
W_{k,n}=\f1{\rho_{k,n}} \sum_{\ell=1}^{\rho_{k,n}} \sum_{\Y\subset \Pnr_{Q_\ell}, \, |\Y|=k} s_n(\Y, \Pnr_{Q_\ell}; t)\, \delta_{H_n(\Y)} \in M_+(E'). 
$$
It then follows from the same argument as Proposition \ref{p:LDP.leading.term} that $(W_{k,n})_{n\ge1}$ satisfies an LDP in the weak topology with rate $\rho_{k,n}$ and the rate function $\bar \Lambda_k^*$. 
Subsequently, by repeating the proof of Proposition \ref{p:exp.negligible}, one can also establish that for every $\delta>0$, 
$$
\lim_{n\to\infty} \f1{\rho_{k,n}} \log \P \big( d_{\ms{TV}} ( V_{k,n}, W_{k,n} )   \ge \delta  \big) = -\infty. 
$$
This concludes the proof of Proposition \ref{p:LDP.Vkn}. 
\end{proof}

\begin{proof}[Proof of Theorem \ref{t:LDP.pp.Ustat}]
{We now  have to show  that $(U_{k,n})_{n\ge1}$ exhibits the same LDP as $(V_{k,n})_{n\ge1}$ above. We take,} without loss of generality, $0<t_1\le t_2 \le \dots \le t_m <\infty$ for time parameters of $(U_{k,n})_{n\ge1}$, whereas we fix  the {parameter of $(V_{k,n})_{n\ge1}$}   at $t=t_1$. In this setup, we need to verity that for every $\delta>0$, 
$$
\lim_{n\to\infty}\f1{\rho_{k,n}} \log \P \big(  { d_{\ms{TV}}(U_{k,n}, V_{k,n})}\ge \delta\big)=-\infty. 
$$
It is straightforward to see that 
$$
{ d_{\ms{TV}}(U_{k,n}, V_{k,n})} \le \f1{\rho_{k,n}}\sum_{\Y\subset \Pn, \, |\Y|=k}s_n(\Y, \Pn; t_1)\big( 1-c_n(\Y, \Pn; t_m) \big). 
$$
By virtue of the partition of $[0,1]^d$ into multiple sub-cubes $Q_1,\dots,Q_{\rho_{k,n}}$ as in the proof of Proposition \ref{p:LDP.leading.term}, what need to be shown are 
\begin{equation}  \label{e:time.parameter.cond1}
\lim_{n\to\infty}\f1{\rho_{k,n}}\log \P \Big( \sum_{\ell=1}^{\rho_{k,n}}\sum_{\Y\subset \Pnr_{Q_\ell}} s_n(\Y, \Pnr_{Q_\ell}; t_1) \big( 1-c_n(\Y, \Pnr_{Q_\ell}; t_m) \big) \ge \delta \rho_{k,n}\Big)=-\infty, 
\end{equation}
and 
\begin{align*}
\lim_{n\to\infty}\f1{\rho_{k,n}}\log \, &\P \Big( \, \Big|  \sum_{\Y\subset \Pn} s_n(\Y, \Pn; t_1)  \big( 1-c_n(\Y, \Pn; t_m) \big) \\
& \qquad - \sum_{\ell=1}^{\rho_{k,n}}\sum_{\Y\subset \Pnr_{Q_\ell}} s_n(\Y, \Pnr_{Q_\ell}; t_1) \big( 1-c_n(\Y, \Pnr_{Q_\ell}; t_m) \big)  \Big| \ge \delta \rho_{k,n} \Big)=-\infty, 
\end{align*}
for every $\delta>0$. Of the last two conditions, the latter  can be established by the same argument as Proposition \ref{p:exp.negligible}. By Markov's inequality, \eqref{e:time.parameter.cond1} will follow if we can show that for every $a>0$, 
\begin{equation}  \label{e:moment.conv.to.1}
\E \Big[ e^{a\sum_{\Y\subset \Pnr_{Q_1}}s_n(\Y, \Pnr_{Q_1}; t_1)(1-c_n(\Y, \Pnr_{Q_1}; t_m))} \Big] \to 1, \ \ n\to\infty. 
\end{equation}
Repeating the same calculations as in \eqref{e:comp.involved}, while using an obvious bound $s_n(\Y, \Pnr_{Q_1}; t_1) \le \one \big\{ \diam(\Y) \le r_nL \big\}$, as well as the Mecke formula for Poisson point processes, 
\begin{align*}
&\E\Big[ \sum_{\Y\subset \Pnr_{Q_1}}s_n(\Y, \Pnr_{Q_1}; t_1)(1-c_n(\Y, \Pnr_{Q_1}; t_m)) \Big]   \\
&\le \f1{k!}\, \int_{(\R^d)^{k-1}} \one \big\{ \diam(\bfz_d,\by)\le L \big\} \Big( 1-e^{-nr_n^d \text{vol} ( \B (\{ \bfz_d,\by \}; t_m) )} \Big)\dif \by \to 0, \ \ \text{as } n\to\infty. 
\end{align*}
Since we have already shown the  uniform integrability by \eqref{e:UI1}, we now obtain \eqref{e:moment.conv.to.1}, as desired. 

Finally, by the same proof as Proposition \ref{p:relative.entropy.and.rate}, we can obtain \eqref{e:rate.function.RE.pp.Ustat}. 
\end{proof}

\subsection{Proofs of  Theorem \ref{t:LDP.Ustat}, Corollary \ref{c:LDP.Ustat}, and Proposition \ref{p:rate.function}} \label{sec:LDP.Ustat.and.rate.function}

In this section, we first prove Proposition \ref{p:rate.function}, because the proof of Theorem \ref{t:LDP.Ustat} makes use of Proposition \ref{p:rate.function}.

\begin{proof}[Proof of Proposition \ref{p:rate.function}]
For every $\bx=(x_1,\dots, x_m)\in E'$, one can choose $\ba=(a_1,\dots, a_m)\in \R^m$ so that 
$$
u_i(\bx, \ba):=\frac{1}{k!}\int_{(\R^d)^{k-1}} h^{(i)} (\bfz_d,\by) e^{\sum_{r=1}^m a_r h^{(r)}(\bfz_d,\by)} \dif \by -x_i = 0, \ \ i=1,\dots,m, 
$$
with $\by\in (\R^d)^{k-1}$. 
By condition \textbf{(H4)}, $u_i$ is continuously differentiable with respect to $(\bx, \ba)$. Furthermore,  \textbf{(H5)} ensures that the matrix $\big( \partial u_i / \partial a_j\big)_{i,j=1}^m$ is positive definite. Thus, by the implicit function theorem, there exist continuously differentiable functions $\alpha_r:E'\to \R$, $r=1,\dots,m$, such that 
\begin{equation}  \label{e:Ik.explicit}
I_k(\bx)=\sum_{i=1}^m \alpha_i(\bx) x_i -\f1{k!}\, \int_{(\R^d)^{k-1}} \Big( e^{\sum_{i=1}^m \alpha_i(\bx)\hi (\bfz_d,\by)}-1 \Big)\dif \by, 
\end{equation}
and 
$$
x_i = \frac{1}{k!}\int_{(\R^d)^{k-1}} h^{(i)} (\bfz_d,\by) e^{\sum_{r=1}^m \alpha_r(\bx) h^{(r)}(\bfz_d,\by)} \dif \by, \ \ i=1,\dots,m. 
$$
In particular, \eqref{e:Ik.explicit} implies that $I_k$ is continuously differentiable on $E'$.

For Part $(ii)$, we show that the Hessian of $I_k$, denoted $\mathcal H(I_k)$, is positive definite. First it is elementary to check that $\partial I_k(\bx)/\partial x_j = \alpha_j(\bx)$ for $j=1,\dots,m$. Using this, it is not hard to calculate that $\mathcal H(I_k) = A^{-1}$ where $A=(A_{ij})$ is given by 
$$
	A_{ij} = \f1{k!}\, \int_{(\R^d)^{k-1}}\hi (\bfz_d,\by) h^{(j)}(\bfz_d,\by) e^{{ \sum_{r=1}^m \alpha_r(\bx)h^{(r)} (\bfz_d,\by)}}\dif \by. 
$$
	We conclude from {condition {\bf (H5)}} that $A$ is positive definite; thus, so is $A^{-1}$ as required. 
Finally, if we set $0=\partial I_k(\bx)/\partial x_j$ for every $j=1,\dots,m$, then $\alpha_j(\bx)=0$, and so, $x_j=\mu_j$ and $I_k(\bx)=0$. As $I_k$ is strictly convex, $(\mu_1,\dots,\mu_m)$ is a unique minimizer of $I_k$. 
\end{proof}

{The proof of Theorem \ref{t:LDP.Ustat} is  based on} an  extension of the contraction principle, provided in \cite[Theorem 4.2.23]{dembo:zeitouni:1998}. We begin with defining a map $F:M_+(E')\to [0,\infty)^m$ by $F(\rho)=\big( \int_{E'} x_i \rho(\dif \bx) \big)_{i=1}^m$. Then, $F(U_{k,n})=T_{k,n}/\rho_{k,n}$, where $U_{k,n}$ and $T_{k,n}$ are  defined in \eqref{e:def.Ukn} and \eqref{e:def.Tkn}, respectively. Since $F$  is not continuous in the weak topology,  {we need to} introduce a family of continuous maps: $F_K(\rho):=\big( \int_{E'} s_i^{(K)}(\bx) \rho(\dif \bx)\big)_{i=1}^m$, where 
$$
s_i^{(K)}(x_1,\dots,x_m) = \begin{cases}
x_i & \text{ if } 0\le x_i \le K, \\
-K^2(x_i-K)+K & \text{ if } K \le x_i \le K+K^{-1}, \\
0 & \text{ if } x_i \ge K+K^{-1}.  
\end{cases}
$$
Clearly, $s_i^{(K)}$ is continuous and bounded on $E'$, and consequently,  $F_K$ becomes  continuous in the weak topology. 

{In order to apply \cite[Theorem 4.2.23]{dembo:zeitouni:1998}, we need to demonstrate the following  auxiliary results}.

\begin{lemma}   \label{l:ustat}
{$(i)$ Let $\delta > 0$. Then,}
	\begin{equation*}
\limsup_{K\to\infty}\limsup_{n\to\infty} \f1{\rho_{k,n}} \log \P \Big( \big\| F(U_{k,n})-F_K(U_{k,n})  \big\| >\delta\Big)=-\infty, 
\end{equation*}
{$(ii)$ Let $a > 0$. Then,}
\begin{equation*}
\limsup_{K\to\infty} \sup_{\rho\in M_+(E'), \, \bar \Lambda_k^*(\rho)\le a} \big\| F(\rho)-F_K(\rho)  \big\|=0. 
\end{equation*}
\end{lemma}
{Assuming the assertions of Lemma \ref{l:ustat} temporarily, we first explain how to conclude the proof of Theorem \ref{t:LDP.Ustat}}.

\begin{proof}[Proof of Theorem \ref{t:LDP.Ustat}]
{Note that one can  conclude from Theorem \ref{t:LDP.pp.Ustat}, \cite[Theorem 4.2.23]{dembo:zeitouni:1998}, and Lemma \ref{l:ustat} that} $(T_{k,n}/\rho_{k,n})_{n\ge1}$ satisfies an LDP with rate $\rho_{k,n}$ and the rate function 
$$
\inf_{\nu\in M_+(E'), \, F(\nu)=\bx}H_k'(\nu|\tau_k'), \ \ \ \bx \in \R^m, 
$$
where $H_k'$ is a relative entropy given at \eqref{e:def.Hk'}.
{It thus remains} to show that 
\begin{equation}  \label{e:rate.func.equal.Ik}
\inf_{\nu\in M_+(E'), \, F(\nu)=\bx}H_k'(\nu|\tau_k')=I_k(\bx), \ \ \ \bx \in \R^m. 
\end{equation}
If $x_i<0$ for some $i\in \{ 1,\dots,m\}$, both sides above are equal to infinity.  So from onward, we consider  $\bx \in E'$ and prove first  that $\inf_{\nu\in M_+(E'), \, F(\nu)=\bx}H_k'(\nu|\tau_k')\ge I_k(\bx)$. 
Fix $\nu\in M_+(E')$ such that $F(\nu)=\bx$. For every $\ba \in \R^m$, we set $f(\bz)=\langle \ba, \bz \rangle$ for  $\bz \in E'$. Approximating $f$ via a sequence of continuous and bounded functions, we have 
\begin{align*}
H_k'(\nu|\tau_k') &= \bar \Lambda_k^*(\nu) \ge \langle \nu, f \rangle -\f1{k!} \int_{(\R^d)^{k-1}} \big(  e^{f(H(\bfz_d,\by))}-1\big)\, \one \big\{ H(\bfz_d,\by)\neq \bfz_m \big\}\dif \by  \\
&= \langle \ba, \bx \rangle -\f1{k!}\int_{(\R^d)^{k-1}} \big( e^{\langle \ba, H(\bfz_d,\by) \rangle}-1 \big)\dif \by; 
\end{align*}
thus, $H_k'(\nu|\tau_k') \ge I_k(\bx)$ holds. 
Next, for every $\bx \in E'$, there exists $\ba_0\in \R^m$ such that 
$$
I_k(\bx) {= \langle \ba_0, \bx \rangle -\f1{k!}\int_{(\R^d)^{k-1}} \big( e^{\langle \ba_0, H(\bfz_d,\by) \rangle}-1 \big)\dif \by} = \langle \ba_0, \bx \rangle -\int_{E'} \big(e^{\langle \ba_0, \bz \rangle}-1  \big)\tau_k'(\dif \bz), 
$$
where $x_i=\int_{E'}z_i e^{\langle \ba_0, \bz \rangle}\tau_k'(\dif \bz)$ for $i=1,\dots,m$. Define $\nu(A):=\int_Ae^{\langle \ba_0, \bz \rangle} \tau_k'(\dif \bz)$, $A\subset E'$. It then follows that $\int_{E'}z_i\nu(\dif \bz)=x_i$ for  $i=1,\dots,m$; equivalently, $F(\nu)=\bx$. Moreover, 
$$
	H_k'(\nu|\tau_k')= \hspace{-.1cm}\int_{E'} \hspace{-.1cm}\log \f{\dif \nu}{\dif \tau_k'} (\bz)\nu (\dif \bz) -\nu(E') + \tau_k'(E')=  \langle \ba_0, \bx \rangle -\hspace{-.1cm}\int_{E'}\hspace{-.1cm} \big( e^{\langle \ba_0, \bz \rangle}-1 \big)\tau_k'(\dif \bz) = I_k(\bx). 
$$
This implies  $\inf_{\nu\in M_+(E'), \, F(\nu)=\bx}H_k'(\nu|\tau_k')\le  I_k(\bx)$, and thus,  the proof of \eqref{e:rate.func.equal.Ik} has been completed. 

Finally, Proposition \ref{p:rate.function} $(i)$ ensures that $I_k$ is continuous; therefore, the LDP for $(T_{k,n}/\rho_{k,n})$ implies the convergence in \eqref{e:LDP.and.continuity}. 

\enp

{Now, we present the proof of Lemma \ref{l:ustat}.}

\bep[Proof of {Lemma \ref{l:ustat} $(i)$}]
	Note first  that the result follows from 
\begin{equation}  \label{e:2nd.prime.cond}
\limsup_{K\to\infty}\limsup_{n\to\infty} \f1{\rho_{k,n}} \log \P \Big( \sum_{\Y\subset \Pn}g_n^{(i)}(\Y, \Pn; t_i)\, \one \big\{ h_n^{(i)}(\Y)>K \big\} \ge \delta \rho_{k,n} \Big)=-\infty, 
\end{equation}
for each $i\in\{ 1,\dots,m \}$. 
In the sequel, we prefer to use the notation $G_n$ and $H_n$, instead of $g_n^{(i)}$ and $h_n^{(i)}$, with the assumption  $m=1$, so that the range of $G_n$ and $H_n$ is $[0,\infty)$. Additionally, we drop the subscript $i$ from $t_i$. 
Then, \eqref{e:2nd.prime.cond} can be rephrased as 
\begin{equation}  \label{e:3rd.prime.cond}
\limsup_{K\to\infty}\limsup_{n\to\infty} \f1{\rho_{k,n}} \log \P \Big( \sum_{\Y\subset \Pn}G_n(\Y, \Pn; t)\, \one \big\{ H_n(\Y)>K \big\} \ge \delta \rho_{k,n} \Big)=-\infty. 
\end{equation}
Clearly, \eqref{e:3rd.prime.cond} is implied by 
\begin{equation}  \label{e:3rd.cond}
\limsup_{K\to\infty}\limsup_{n\to\infty} \f1{\rho_{k,n}} \log \P \Big( \sum_{\ell=1}^{\rho_{k,n}} \sum_{\Y\subset \Pnr_{Q_\ell}} G_n(\Y, \Pnr_{Q_\ell}; t)\, \one \big\{  H_n(\Y)>K \big\} \ge \delta \rho_{k,n}\Big)=-\infty, 
\end{equation}
and 
\begin{align}  
&\limsup_{K\to\infty}\limsup_{n\to\infty} \f1{\rho_{k,n}} \log \P \bigg( \, \Big| \sum_{\Y\subset \Pn} G_n(\Y, \Pn; t)\, \one \big\{ H_n(\Y)>K \big\}  \label{e:4th.cond} \\
&\qquad \qquad \qquad \qquad \qquad  -\sum_{\ell=1}^{\rho_{k,n}} \sum_{\Y\subset \Pnr_{Q_\ell}} G_n(\Y, \Pnr_{Q_\ell}; t)\, \one \big\{  H_n(\Y)>K \big\}  \Big| \ge \delta \rho_{k,n}  \bigg)=-\infty,  \notag
\end{align}
for every $\delta>0$. As for \eqref{e:3rd.cond}, Markov's inequality yields that 
\begin{align*}
&\f1{\rho_{k,n}} \log \P \Big( \sum_{\ell=1}^{\rho_{k,n}} \sum_{\Y\subset \Pnr_{Q_\ell}} G_n(\Y, \Pnr_{Q_\ell}; t)\, \one \big\{  H_n(\Y)>K \big\} \ge \delta \rho_{k,n}\Big)  \\
&\le -a\delta + \log \E \Big[ e^{  a\sum_{\Y\subset \Pnr_{Q_1}} G_n (\Y, \Pnr_{Q_1}; t)\, \one \{ H_n(\Y)>K \} } \Big]. 
\end{align*}
{
It is thus sufficient to demonstrate that for every $a>0$, 
\begin{equation}  \label{e:UI.and.K}
\limsup_{K\to\infty} \limsup_{n\to\infty} \E \Big[  e^{a\sum_{\Y\subset \Pnr_{Q_1}} G_n(\Y, \Pnr_{Q_1}; t) \one \{ H_n(\Y)>K \}} \Big] = 1. 
\end{equation}
Although the proof techniques for \eqref{e:UI.and.K} are mostly the same  as those for \eqref{e:UI1},  
we still need  to update it because $G_n$ in \eqref{e:UI.and.K}  is not necessarily  an indicator function. Using the same logic as in the proof of \eqref{e:UI1}, our goal is to prove that for every $a>0$, 
\begin{equation*}  \label{e:suff.cond.UI.non.indicator}
\limsup_{K\to\infty}\limsup_{n\to\infty} \bigg( \E\Big[ e^{a\sum_{\Y\subset \Pnr_{Q_1}} G_n(\Y, \Pnr_{Q_1}; t) \one \{ H_n(\Y)>K, \, \ell(\Y)\in J_1 \}} \Big] \bigg)^{b_n} \le 1, 
\end{equation*}
where $J_1=\big[ 0,tr_n/\sqrt{d} \big]^d$ and $b_n=\rho_{k,n}^{-1}/(4Lr_n)^d$. Observe, once again, that there exists at most a single $k$-point set $\Y\subset \Pnr_{Q_1}$ satisfying $c_n(\Y, \Pnr_{Q_1}; t)=1$ and $\ell(\Y)\in J_1$.  Therefore, as $n\to\infty$, 
\begin{align}
&\E\Big[ e^{a\sum_{\Y\subset \Pnr_{Q_1}} G_n(\Y, \Pnr_{Q_1}; t) \one \{H_n(\Y)>K, \, \ell(\Y)\in J_1 \}} \Big]  \label{e:UI.K.expectation}\\
&\le 1 + \E \Big[ \sum_{\Y\subset \Pn} e^{aH_n(\Y)} \one \big\{ H_n(\Y) >K, \, \ell(\Y)\in J_1 \big\} \Big] \notag \\
&=1+\f{n^k}{k!}\int_{([0,1]^d)^k} e^{aH_n(x_1,\dots,x_k)} \one \big\{ H_n(x_1,\dots,x_k)>K, \, \ell(x_1,\dots,x_k)\in J_1 \big\}\dif \bx \notag \\
&=1+\f{\rho_{k,n}}{k!}\int_{[0,1]^d}\int_{(\R^d)^{k-1}} e^{aH(\bfz_d,\by)} \one \big\{  H(\bfz_d,\by)>K, \, \ell(x,x+r_ny_1,\dots,x+r_ny_{k-1})\in J_1\big\} \notag \\
&\qquad \qquad \qquad \qquad \qquad \qquad \qquad \qquad \qquad \qquad  \times \prod_{i=1}^{k-1}\one \big\{ x+r_ny_i\in [0,1]^d \big\}\dif \by \dif x \notag \\
&\le 1+C^* \rho_{k,n} r_n^d \int_{(\R^d)^{k-1}} e^{aH(\bfz_d,\by)}\one \big\{  H(\bfz_d,\by)>K\big\}\dif \by. \notag 
\end{align}
Here, the last integral is finite due to property \textbf{(H4)}. Now, it follows from \eqref{e:UI.K.expectation} that 
\begin{align*}
&\limsup_{K\to\infty}\limsup_{n\to\infty} \bigg( \E\Big[ e^{a\sum_{\Y\subset \Pnr_{Q_1}} G_n(\Y, \Pnr_{Q_1}; t) \one \{ H_n(\Y)>K, \, \ell(\Y)\in J_1 \}} \Big] \bigg)^{b_n} \\
&\le \limsup_{K\to\infty} \exp \Big\{ \frac{C^*}{(4L)^d} \int_{(\R^d)^{k-1}} e^{aH(\bfz_d,\by)}\one \big\{ H_n(\Y) >K\big\} \dif \by \Big\} =1
\end{align*}
}

Next, turning to condition \eqref{e:4th.cond}, we deduce an  inequality analogous to \eqref{e:disect.approx.two}:
\begin{align*}
&\Big| \sum_{\Y\subset \Pn} G_n(\Y, \Pn; t)\, \one \big\{ H_n(\Y)>K \big\} -\sum_{\ell=1}^{\rho_{k,n}} \sum_{\Y\subset \Pnr_{Q_\ell}} G_n(\Y, \Pnr_{Q_\ell}; t)\, \one \big\{  H_n(\Y)>K \big\}  \Big| \\
&\le \sum_{\Y\subset \Pn} G_n \big(\Y, \Pnr_{\bigcup_{\ell=1}^{\rho_{k,n}} Q_{\ell}^\partial ((L+t)r_n)}; t \big) \one \Big\{ \Y\subset \bigcup_{\ell=1}^{\rho_{k,n}} Q_\ell^\partial \big((L+t)r_n\big) \Big\}    \\
&\qquad \qquad \qquad \qquad \qquad \qquad  + \sum_{\ell=1}^{\rho_{k,n}} \sum_{\Y\subset \Pnr_{\bQ_\ell((L+t)r_n)} } G_n \big(\Y, \Pnr_{\bQ_\ell ((L+t)r_n)}; t  \big).
\end{align*}
Proceeding the same argument as those after \eqref{e:disect.approx.two}, while using property \textbf{(H4)}, we can get \eqref{e:4th.cond}. Now, the proof of {Lemma \ref{l:ustat} $(i)$} is completed. 
\enp

\bep[Proof of {Lemma \ref{l:ustat} $(ii)$}]
Suppose, for contradiction, that we can choose  $a>0$ and $\delta>0$ such that 
\begin{equation}  \label{e:2nd.cond.contradiction}
\limsup_{K\to\infty} \sup_{\rho\in M_+(E'), \, \bar \Lambda_k^*(\rho)\le a} \big\| F(\rho)-F_K(\rho)  \big\| >\delta. 
\end{equation}
For $K'\ge K+K^{-1}$, define
$$
F_{K, K'}(\rho):=\Big( \int_{E'}u_i^{(K, K')}(\bx)\rho(\dif \bx) \Big)_{i=1}^m, \  \ \ \rho\in M_+(E'), 
$$
where 
$$
u_i^{(K, K')}(x_1,\dots,x_m) = \begin{cases}
0 & \text{ if } x_i \in [0,K] \cup [K'+(K')^{-1},\infty), \\
(K^2+1)(x_i-K) & \text{ if } K \le x_i \le K+K^{-1}, \\
x_i & \text{ if } K+K^{-1} \le x_i \le K', \\
-(K')^2 (x_i-K')+K' & \text{ if } K' \le x_i \le K' + (K')^{-1}. 
\end{cases}
$$
Then, $F_{K, K'}$ is continuous in the weak topology, and further, 
\begin{equation}  \label{e:FKK'.conv}
\big\| F_{K, K'}(\rho) \big\| \nearrow \big\| F(\rho)-F_K(\rho) \big\|, \ \ \text{as } K'\uparrow \infty. 
\end{equation}
{By Lemma \ref{l:ustat} $(i)$} and \eqref{e:FKK'.conv}, we can choose $K_1>0$ so that 
\begin{equation}  \label{e:LDP.FKK'}
\limsup_{n\to\infty} \f1{\rho_{k,n}} \log \P \Big( \big\| F_{K, K'}(U_{k,n}) \big\| >\delta \Big) <-a
\end{equation}
for any $K'$ and $K$, with $K'\ge K+K^{-1} > K \ge K_1$. By virtue of the continuity of $F_{K, K'}$ together with \eqref{e:LDP.FKK'}, one can apply the contraction principle to the LDP for $(U_{k,n})_{n\ge1}$ (see Theorem \ref{t:LDP.pp.Ustat}), to obtain that 
$$
\inf_{\rho\in M_+(E'), \, \| F_{K, K'}(\rho) \|>\delta}\bar \Lambda_k^*(\rho)>a, 
$$
for all $K'\ge K+K^{-1} > K \ge K_1$. 

Next, let us turn to \eqref{e:2nd.cond.contradiction} and  fix $K\ge K_1$, so that $\sup_{\rho\in M_+(E'), \, \bar \Lambda_k^*(\rho)\le a}\big\| F(\rho)-F_K(\rho) \big\|>\delta$. Then, we can choose $\nu\in M_+(E')$ with $\bar \Lambda_k^*(\nu)\le a$ and $\big\| F(\nu)-F_K(\nu) \big\|>\delta$. Furthermore, \eqref{e:FKK'.conv} implies that there exists $K'\ge K+K^{-1}$ with $\big\| F_{K, K'}(\nu) \big\|>\delta$. Since this is contradiction, we have established {Lemma \ref{l:ustat} $(ii)$}, as desired. 
\end{proof}

{Finally, we prove Corollary \ref{c:LDP.Ustat}. The proof  is very similar to that of Corollary \ref{c:LDP.pp.general}. Hence, we  focus only on some of its key differences. The goal is to show that for every $\varepsilon_0 > 0$,}  
$$\lim_{n \to \infty}\rho_{k, n}^{-1}\log\P\big(\|T_{k, n} - T_{k,n}^{\ms B}\|\ge \varepsilon_0\rho_{k,n}\big) = -\infty.$$

\begin{proof}[Proof of Corollary \ref{c:LDP.Ustat}]
As in the proof of Lemma \ref{l:ustat} $(i)$, we may put the assumption $m=1$ and use the notation $G_n$ and $H_n$ (instead of $g_n^{(1)}$ and $h_n^{(1)}$). 
For $K > 0$, we decompose  $T_{k,n}^{\ms B} = T_{\le K, k, n}^{\ms B} + T_{> K, k, n}^{\ms B}$, where
$$	T_{\le K,k,n}^{\ms B}:= \sum_{\Y\subset \Bn, \, |\Y|=k} G_n(\Y, \Bn; t) \one \{H_n(\Y) \le K\},	$$
	and 
$$	T_{> K,k,n}^{\ms B}:= \sum_{\Y\subset \Bn, \, |\Y|=k} G_n(\Y, \Bn; t) \one \{H_n(\Y) > K\}.	$$
Similarly, we write $T_{k, n} = T_{\le K, k, n} + T_{> K, k, n}$. 
Since all the summands in $T_{\le K,k,n}$ and $T_{\le K,k,n}^{\ms B}$ are  bounded by $K$, a simple repetition of the proof of Corollary \ref{c:LDP.pp.general} can yield that 
$$
\frac{1}{\rho_{k,n}} \log \P \big( \, |\, T_{\le K,k,n} - T_{\le K,k,n}^{\ms B} \,| \ge \varepsilon_0 \rho_{k,n} \big) \to -\infty, \ \ \ n\to\infty, 
$$
for every $K>0$. Additionally, the proof of Lemma \ref{l:ustat} $(i)$ has already shown that 
$$
\limsup_{K\to\infty} \limsup_{n\to\infty} \frac{1}{\rho_{k,n}} \log \P ( T_{>K,k,n}  \ge \varepsilon_0 \rho_{k,n}) = -\infty; 
$$
so it  remains to verify that 
\begin{equation}  \label{e:negligible.part.B}
\limsup_{K \to \infty}\limsup_{n \to \infty}\frac{1}{\rho_{k,n}} \log\P\big(T_{> K, k, n}^{\ms B} \ge \varepsilon_0 \rho_{k,n}\big) = -\infty.
\end{equation}
Using the same diluted cubes  from the proof of Corollary \ref{c:LDP.pp.general}, \eqref{e:negligible.part.B} will follow, provided that 
\begin{align*}
&\limsup_{K\to\infty} \limsup_{n\to\infty} \frac{1}{\rho_{k,n}} \log \P  \Big( \sum_{i=1}^{b_n'}\sum_{\Y \subset \Bn, \, |\Y|=k} G_n(\Y, \Bn; t)\\
&\qquad\qquad \qquad \qquad \qquad \qquad\qquad \times  \one \big\{  H_n(\Y) >K, \, \ell(\Y) \in J_i \big\} \ge \varepsilon_0 \rho_{k,n} \Big) = -\infty. 
\end{align*}
 Observe again that there exists a  constant $M > 0$ such that 
	\begin{align}  \label{e:max.bound}
	\begin{split}
		&\sum_{\Y\subset \Bn, \, |\Y|=k} G_n(\Y, \Bn; t) \one \{H_n(\Y) > K, \, \ell(\Y) \in J_i\}\\
		&\quad \le M \max_{\Y\subset \Bn, \, |\Y|=k}\Big( H_n(\Y) \one \{H_n(\Y) > K, \, \ell(\Y) \in J_i\}\Big).
		\end{split}
	\end{align}
	Moreover, as in the proof of Corollary \ref{c:LDP.pp.general}, we may work under the event  $\big\{\Bn \subset \Pn^{(1,\ms{a})}\big\}$, where $\Pn^{(1,\ms{a})}=\Pn \cup \Pn^{(1)}\stackrel{d}{=}\mathcal P_{2n}$. 
	By Markov's inequality and \eqref{e:max.bound}, we have,   for  every $a > 0$,
	\begin{align*}
		&\P\Big(\Big\{ \sum_{i=1}^{b_n'}\sum_{\Y\subset \Bn, \, |\Y|=k} \hspace{-5pt}G_n(\Y, \Bn; t) \one \{H_n(\Y) > K, \, \ell(\Y) \in J_i\} \ge \varepsilon_0\rho_{k, n} \Big\} \cap \big\{\Bn \subset \Pn^{(1,\ms{a})}\big\} \Big) \\
		&\quad \le \P \Big( \sum_{i=1}^{b_n'}  M \max_{\Y\subset \Pn^{(1,\ms{a})}, \, |\Y|=k}\Big( H_n(\Y) \one \{H_n(\Y) > K, \, \ell(\Y) \in J_i\}\Big) \ge \varepsilon_0\rho_{k, n} \Big) \\
		&\quad\le e^{-a\varepsilon_0\rho_{k, n}} \big(\E\big[e^{aZ_n}\big]\big)^{b_n'}, 
	\end{align*}
	where 
	$$Z_n :=M \hspace{-.4cm}\max_{\Y\subset \mathcal P_{2n}, \, |\Y|=k}\Big( H_n(\Y) \one \{H_n(\Y) > K, \, \ell(\Y) \in J_1\}\Big).$$
Now, proceeding as in \eqref{e:UI.K.expectation}, 
	\begin{align*}
	\E\big[e^{aZ_n}\big]
&\le 1 + \E \Big[ \sum_{\Y\subset \mathcal P_{2n}, \, |\Y|=k} e^{aMH_n(\Y)} \one \big\{ H_n(\Y) >K, \, \ell(\Y)\in J_1 \big\} \Big] \\
		&\le 1+ C^* \rho_{k,n} r_n^d \int_{(\R^d)^{k-1}}\hspace{-.5cm} e^{aMH(\bfz_d,\by)}\one \big\{  H(\bfz_d,\by)>K\big\}\dif \by \\
		&\le \exp \Big\{ C^* \rho_{k,n} r_n^d \int_{(\R^d)^{k-1}}\hspace{-.5cm} e^{aMH(\bfz_d,\by)}\one \big\{  H(\bfz_d,\by)>K\big\}\dif \by  \Big\}. 
	\end{align*}
{This, together with property \textbf{(H4)}, implies} that
	\begin{align*}
&\limsup_{K\to \infty}	\limsup_{n\to\infty}\frac{1}{\rho_{k, n}}	\log\Big\{ e^{-a\varepsilon_0\rho_{k, n}} \Big(\E\big[e^{aZ_n}\big]\Big)^{b_n'}\Big\}  \\
&\quad \le -a\varepsilon_0 + C^*\limsup_{K\to \infty}	\int_{(\R^d)^{k-1}}\hspace{-.5cm} e^{aMH(\bfz_d,\by)}\one \big\{  H(\bfz_d,\by)>K\big\}\dif \by = -a\varepsilon_0.
	\end{align*}
	Since $a$ is arbitrary, letting $a\to\infty$ concludes the entire proof. 
\end{proof}
\subsection{Proofs of Theorem \ref{t:LDP.persistent.Betti} and  {Theorem} \ref{t:LDP.critical.point}}  \label{sec:LDP.critical.point.and.LDP.persistent.Betti}

\begin{proof}[Proof of Theorem \ref{t:LDP.persistent.Betti} (Poisson input)] {We first deal with the case of a Poisson input.}
Since the function $\big( h_{s_i}(x_1,x_2,x_3)h_{t_i}(x_1,x_2,x_3) \big)_{i=1}^m$ defined at \eqref{e:def.h.persistent.betti} satisfies {conditions \textbf{(H1)}--\textbf{(H5)}}, 
a direct application of Theorem  \ref{t:LDP.Ustat} yields that as $n\to\infty$, 
$$
\f1{\rho_{3,n}}\, \log \P \bigg(  \Big(  \rho_{3,n}^{-1} \sum_{\Y\subset \Pn, \, |\Y|=3} g_{r_ns_i, r_nt_i} (\Y, \Pn), \, i=1,\dots,m \Big)\in A \bigg) \to -\inf_{\bx\in A}I_3(\bx). 
$$
To complete the proof we show that for every $0 \le s \le t <\infty$ and $\delta>0$, 
$$
\lim_{n\to\infty} \f1{\rho_{3,n}}\, \log \P \bigg(  \beta_{1,n}(s,t) - \sum_{\Y\subset \Pn, \, |\Y|=3} g_{r_ns, r_nt} (\Y, \Pn) \ge \delta \rho_{3,n} \bigg)=-\infty. 
$$
By definition, $\beta_{1,n}(s,t)$ represents the number of persistent $1$-cycles in the region $[0,s]\times[t,\infty]$ of the first-order persistence diagram. In particular, $\beta_{1,n}(s,t)$ accounts for subsets of $3$ points in $\R^2$ that form a single $1$-cycle before time $r_ns$, such that this 1-cycle remains alive at time $r_nt$ and isolated from all the remaining points in $\Pn$ at that time. Note that these subsets  of $3$ points in $\R^2$ are counted  by $\sum_{\Y\subset \Pn, \, |\Y|=3} g_{r_ns, r_nt} (\Y, \Pn)$ as well. Thus,  all the remaining points in $[0,s]\times [t,\infty]$ of the first-order persistence diagram are  associated to the 1-cycles  on connected components of size greater than or equal to $4$ at time $r_nt$. From this point of view, 
\begin{equation*}  \label{e:difference.persistent.Betti.leading}
\beta_{1,n}(s,t) - \sum_{\Y\subset \Pn, \, |\Y|=3} g_{r_ns, r_nt} (\Y, \Pn)
\end{equation*}
can be bounded by the first-order Betti number at time $r_nt$, associated only to connected components of size greater than or equal to $4$. Moreover, this Betti number  is  further  bounded by the corresponding $1$-simplex counts {(i.e., edge counts). By Lemma 2.1 in \cite{dereudre:georgii:2009}, one can bound such $1$-simplex counts}  by three times the number of vertices that are contained in connected components of size greater than or equal to $4$. In conclusion,  we have that
\begin{equation}  \label{e:d=2.k=3.Betti}
\beta_{1,n}(s,t) - \sum_{\Y\subset \Pn, \, |\Y|=3} g_{r_ns, r_nt} (\Y, \Pn)\le 3 \sum_{i=4}^\infty i \sum_{\Y\subset \Pn, \, |\Y|=i} \vni (\Y, \Pn), 
\end{equation}
where 
$$
\vni (\Y, \Pn)=\one \big\{ \alpha(\Y, r_nt) \text{ is a connected component of } \alpha(\Pn, r_nt) \big\}. 
$$

As in the proof of Proposition \ref{p:LDP.leading.term}, we partition $[0,1]^2$ into sub-cubes $Q_1,\dots,Q_{\rho_{3,n}}$ of volume $\rho_{3,n}^{-1}$ so that  $Q_1=[0,\rho_{3,n}^{-1/2}]^2$. {We claim that 
\begin{equation}  \label{e:further.new.bound}
\sum_{i=4}^\infty i \sum_{\Y\subset \Pn, \, |\Y|=i} \hspace{-5pt}v_n^{(i)} (\Y, \Pn) \le \sum_{\bz\in \{ 0,\pm 1 \}^2}\sum_{\ell=1}^{\rho_{3,n}} \sum_{i=4}^\infty i \sum_{\Y \subset \Pnr_{(Q_\ell + 3r_nt\bz)\cap [0,1]^2}} \hspace{-20pt}v_n^{(i)} \big( \Y, \Pnr_{(Q_\ell + 3r_nt\bz)\cap [0,1]^2} \big). 
\end{equation}
To show \eqref{e:further.new.bound}, suppose there exists a connected component $\Y$ of cardinality at least $4$, whose points are  counted by the left hand side of \eqref{e:further.new.bound}. If $\Y$ is contained in one of the cubes in $(Q_\ell)_{\ell=1}^{\rho_{3,n}}$, all the points in $\Y$ can also  be counted by the statistics on the right hand side of \eqref{e:further.new.bound} with $\bz=(0,0)$. If we observe a subset $\Z\subset \Y \cap Q_\ell$ for some $\ell$ with $|\Z|\le 3$, such that $\Z$ itself forms a connected component within $Q_\ell$, with respect to the  process $\Pnr_{Q_\ell}$, then the points in $\Z$ will be missed from the above statistics with $\bz=(0,0)$. Even in that case, however, all the points in $\Z$ can eventually be counted by the statistics in \eqref{e:further.new.bound} with other choice of $\bz\in \{  0,\pm1\}^2\setminus ( 0,0 )$. 
Since there are at most finitely many choices of $\bz$, the entire proof will be complete, provided that for any $\delta>0$, 
\begin{equation}  \label{e:approx.persistent.Betti1}
\lim_{n\to\infty}\f1{\rho_{3,n}}\, \log \P \bigg( \sum_{\ell=1}^{\rho_{3,n}} \sum_{i=4}^\infty i \sum_{\Y\subset \Pnr_{Q_\ell}, \, |\Y|=i} \vni \big( \Y, \Pnr_{Q_\ell} \big) \ge \delta \rho_{3,n} \bigg)=-\infty. 
\end{equation}
}
For the proof of  \eqref{e:approx.persistent.Betti1} we apply Markov's inequality to obtain that 
\begin{align*}
&\f1{\rho_{3,n}}\, \log \P \bigg( \sum_{\ell=1}^{\rho_{3,n}} \sum_{i=4}^\infty i \sum_{\Y\subset \Pnr_{Q_\ell}, \, |\Y|=i} \vni \big( \Y, \Pnr_{Q_\ell} \big) \ge \delta \rho_{3,n} \bigg) \\
&\le -a\delta + \log \E \bigg[ \exp \Big\{  a\sum_{i=4}^\infty i \sum_{\Y\subset \Pnr_{Q_1}, \, |\Y|=i} \vni \big( \Y, \Pnr_{Q_1} \big) \Big\}  \bigg]. 
\end{align*}
From this, it suffices to show that for every $a>0$, 
\begin{equation}  \label{e:exp.vanish}
\E \bigg[ \exp \Big\{  a\sum_{i=4}^\infty i \sum_{\Y\subset \Pnr_{Q_1}, \, |\Y|=i} \vni \big( \Y, \Pnr_{Q_1} \big) \Big\}  \bigg] \to 1, \ \ \text{as } n\to\infty. 
\end{equation}
To show this, we first claim that 
\begin{equation}  \label{e:conv.to.zero.in.prob}
\sum_{i=4}^\infty i \sum_{\Y\subset \Pnr_{Q_1}, \, |\Y|=i} \vni \big( \Y, \Pnr_{Q_1} \big) \stackrel{\P}{\to}0, \ \ \ n\to\infty. 
\end{equation}
Taking an expectation of \eqref{e:conv.to.zero.in.prob}, 
\begin{align*}
\sum_{i=4}^\infty i\,  \E \Big[ \sum_{\Y\subset \Pnr_{Q_1}, \, |\Y|=i} \vni \big( \Y, \Pnr_{Q_1} \big)\Big] \le \sum_{i=4}^\infty i\,  \E \Big[ \sum_{\Y\subset \Pnr_{Q_1}, \, |\Y|=i} \uni ( \Y)\Big], 
\end{align*}
where 
$
\uni(\Y): = \one \big\{ \alpha (\Y, r_nt) \text{ is connected} \big\}. $
By the Mecke formula for Poisson point processes, along with the customary change of variables, the right hand side of the above is equal to 
\begin{align*}
&\sum_{i=4}^\infty i\, \f{n^i}{i!}\, \int_{(Q_1)^i} \uni(x_1,\dots,x_i)\dif \bx \\
	&=\sum_{i=4}^\infty  \f{n^ir_n^{2(i-1)}}{(i-1)!}\hspace{-5pt} \int_{Q_1}\hspace{-2pt}\int_{(\R^2)^{i-1}} \hspace{-20pt} \one \big\{ \alpha\big( \{ \bfz_2,y_1,\dots,y_{i-1} \}, t \big) \text{ is connected}  \big\}  \prod_{\ell=1}^{i-1}\one \{ x+r_ny_\ell \in Q_1 \}\dif \by \dif x \\
&\le \sum_{i=4}^\infty \f{(n r_n^2)^{i-3}}{(i-1)!}\, \int_{(\R^2)^{i-1}}\hspace{-5pt}\one \Big\{ \alpha\big( \{ \bfz_2,y_1,\dots,y_{i-1} \}, t \big) \text{ is connected}  \Big\}\dif \by. 
\end{align*}
Because of an elementary fact that there exist at most $i^{i-2}$ spanning trees on $i$ vertices, 
$$
\int_{(\R^2)^{i-1}}\hspace{-5pt}\one \Big\{ \alpha\big( \{ \bfz_2,y_1,\dots,y_{i-1} \}, t \big) \text{ is connected}  \Big\}\dif \by \le i^{i-2}(t^2\pi)^{i-1}. 
$$
 Referring this bound back into the above and appealing to Stirling's formula, i.e., $(i-1)! \ge \big((i-1)/e \big)^{i-1}$ for large $i$, we conclude that as $n\to\infty$, 
\begin{align*}
\sum_{i=4}^\infty i\,  \E \Big[ \sum_{\Y\subset \Pnr_{Q_1}, \, |\Y|=i} \uni ( \Y)\Big] &\le \sum_{i=4}^\infty \f{(n r_n^2)^{i-3}}{(i-1)!} \, i^{i-2} (t^2\pi)^{i-1} \le C^* \sum_{i=4}^\infty \big( et^2 \pi nr_n^2 \big)^{i-3} \to 0. 
\end{align*}
To conclude \eqref{e:exp.vanish} from \eqref{e:conv.to.zero.in.prob}, one needs to show the uniform integrability: for every $a>0$, 
\begin{equation} \label{e:UI.approx}
\limsup_{n\to\infty} \E\bigg[ \exp \Big\{  a\sum_{i=4}^\infty i \sum_{\Y\subset \Pnr_{Q_1}, \, |\Y|=i} \vni \big( \Y, \Pnr_{Q_1} \big) \Big\}  \bigg] <\infty. 
\end{equation}
For the proof, we consider the family $G$ of diluted cubes of side length $8r_n$ that are  contained in $Q_1=\big[ 0,\rho_{3,n}^{-1/2} \big]^2$. Then, the total number of such cubes in $G$ is $c_n:=\rho_{3,n}^{-1}/(8r_n)^2$, which is assumed without loss of generality to be integer-valued for every $n$.  Observe now that for any connected component of size greater than or equal to $4$, there exist $\bz=(z_1, z_2)\in \big\{ 0,\pm1, \dots, \pm4\big\}^2$ and $J\in G$ such that $J+r_n\bz\subset Q_1$ and $\big| \Pn\cap (J+r_n\bz)  \big|\ge 4$. In conclusion, we have 
\begin{align*}
&\sum_{i=4}^\infty i \sum_{\Y\subset \Pnr_{Q_1}, \, |\Y|=i} \vni \big( \Y, \Pnr_{Q_1} \big) \\
&\le \sum_{\bz\in \{ 0,\pm1, \dots, \pm4\}^2} \sum_{J\in G, \, J+r_n\bz\subset Q_1} \big| \Pn\cap (J+r_n\bz) \big|\, \one \Big\{ \big| \Pn\cap (J+r_n\bz) \big| \ge 4 \Big\}. 
\end{align*}
Now, according to H\"older's inequality as well as the homogeneity of $\Pn$, \eqref{e:UI.approx} follows from  
$$
\limsup_{n\to\infty} \E \Big[  e^{a \sum_{J\in G} |\Pn\cap J|\, \one \{ |\Pn\cap J|\ge 4 \} }  \Big]<\infty,
$$
for every $a>0$. Writing $G=\{ J_1,\dots,J_{c_n} \}$ with $J_1=[ 0,8r_n ]^2$, we have that 
$$
\E \Big[  e^{a \sum_{J\in G} |\Pn\cap J|\, \one \{ |\Pn\cap J|\ge 4 \} }  \Big]= \bigg\{ \Big( \E\Big[ e^{a\,|\Pn\cap J_1|\, \one \{ |\Pn\cap J_1|\ge 4 \}} \Big] \Big)^{1/(\rho_{3,n}r_n^2)}  \bigg\}^{1/64}
$$
It is elementary to calculate that 
$$
\E\Big[ e^{a\,|\Pn\cap J_1|\, \one \{ |\Pn\cap J_1|\ge 4 \}} \Big] \le 1 + \sum_{\ell=4}^\infty e^{a\ell} \P\big( |\Pn\cap J_1|=\ell \big) \le 1+ C^* (nr_n^2)^{4}. 
$$
Since $\rho_{3,n}r_n^2=(nr_n^2)^3$, one can obtain that 
\begin{align*}
&\limsup_{n\to\infty} \Big( \E\Big[ e^{t\,|\Pn\cap J_1|\, \one \{ |\Pn\cap J_1|\ge 4 \}} \Big] \Big)^{1/(\rho_{3,n}r_n^2)} \le \limsup_{n\to\infty} \big( 1 + C^*(nr_n^2)^{4}\big)^{1/(nr_n^2)^3} = 1.
\end{align*}

\end{proof}
{
\begin{proof}[Proof of Theorem \ref{t:LDP.persistent.Betti} (binomial input)] 
In the case of a binomial input, instead of \eqref{e:d=2.k=3.Betti}, we  deduce that
$$
\beta_{1,n}^\ms{B}(s,t) - \sum_{\Y\subset \Bn, \, |\Y|=3} g_{r_ns, r_nt} (\Y, \Bn)\le 3 \sum_{i=4}^n i \sum_{\Y\subset \Bn, \, |\Y|=i} \vni (\Y, \Bn).  
$$
As in the proofs of Corollaries \ref{c:LDP.pp.general} and \ref{c:LDP.Ustat}, we may work under the event $\{ \Bn\subset \Pn^{(1,\ms{a})} \}$, where $\Pn^{(1,\ms{a})} = \Pn \cup \Pn^{(1)} \stackrel{d}{=} \mathcal P_{2n}$. Then, for every $\delta>0$, 
\begin{align}
\begin{split}  \label{e:binomial.betti.to.Poisson}
&\frac{1}{\rho_{3,n}} \log \P \bigg( \Big\{  \sum_{i=4}^n i \sum_{\Y\subset \Bn, \, |\Y|=i} \vni (\Y, \Bn) \ge \delta \rho_{3,n}\Big\}  \cap \{ \Bn\subset \Pn^{(1,\ms{a})} \} \bigg) \\
&\le \frac{1}{\rho_{3,n}} \log \P  \Big(  \sum_{i=4}^\infty i \sum_{\Y\subset \mathcal P_{2n}, \, |\Y|=i} \vni (\Y, \mathcal P_{2n}) \ge \delta \rho_{3,n} \Big). 
\end{split}
\end{align}
Here, we used that the double sum counts the number of vertices in connected components of size at least 4. Moreover,  adding further points in $\Pn^{(1,\ms{a})}\setminus \Bn$ increases the number of points in the components associated to $\Bn$. This is clear for the \v Cech complex and follows by the nerve lemma for the alpha complex.
Repeating the same argument as in the Poisson input case, one can show that the right-hand side in \eqref{e:binomial.betti.to.Poisson} goes to $-\infty$ as $n\to\infty$. 
\end{proof}
}

\begin{proof}[Proof of Theorem \ref{t:LDP.critical.point}]
We start by formulating the function $H:=(h^{(1)}, \dots, h^{(m)}): (\R^2)^3 \to [0,\infty)^m$ in the current setup: 
$$
h^{(i)}(x_1, x_2,x_3) :=\one \big\{ \gamma(x_1, x_2,x_3)\in \text{conv}^\circ (x_1, x_2,x_3), \, \mathcal R(x_1, x_2,x_3)\le t_i \big\}. 
$$
Clearly, $H$ satisfies {conditions \textbf{(H1)}--\textbf{(H5)}}. 
In particular,  fix a constant $L$ determined by  property \textbf{(H3)}. 
Define the scaled version $H_n$ of $H$ as in \eqref{e:def.Hn}. For a $3$-point subset $\Y\subset \R^2$ and a finite subset $\Z\supset \Y$ in $\R^2$, define 
\begin{equation}  \label{e:def.c.critical.point}
c(\Y, \Z) := \big( \one \big\{ \U(\Y) \cap \Z = \emptyset \big\} \big)_{i=1}^m. 
\end{equation}
In contrast to \eqref{e:def.c}, the function \eqref{e:def.c.critical.point}  does not involve a time parameter $\bt=(t_1,\dots,t_m)$. Furthermore, unlike \eqref{e:def.cn.and.Gn}, $c_n(\Y,  \Z):=c(r_n^{-1}\Y; r_n^{-1}\Z)=c(\Y, \Z)$ does not depend on $n\ge1$. Finally we  define  $s_n(\Y, \Z):= c_n(\Y, \Z)\, \one \big\{ \diam(\Y)\le r_nL \big\}$ as in \eqref{e:def.sn}. 

The required large deviations in {Theorem \ref{t:LDP.critical.point}} can be deduced by an application of Theorem \ref{t:LDP.Ustat}  to the number of Morse  critical points in \eqref{e:def.Nkn}. Before doing so, however, one must properly modify the proof of Theorem \ref{t:LDP.pp.general} under the setup of {Theorem \ref{t:LDP.critical.point}}. After that, one must also modify the proofs of Theorems \ref{t:LDP.pp.Ustat} and \ref{t:LDP.Ustat}; however, the required modification for these two theorems seems to be sufficiently simple, so we focus our attention only to Theorem \ref{t:LDP.pp.general}. In the below, we  discuss two specific points. 
First,  one has to modify the calculation in \eqref{e:comp.involved} as follows: as $n\to\infty$, 
\begin{align*}
&\f{n^3}{6}\int_{(Q_\ell)^3} \one \big\{ \diam(x_1, x_2,x_3)\le r_nL\big\} \Big( 1-e^{-n\text{vol} \big(\U(x_1, x_2,x_3) \cap Q_\ell  \big)} \Big) \dif \bx  \\
&= \f{\rho_{3,n}}{6}\int_{Q_\ell}\int_{(\R^2)^2} \one \big\{ \diam(\bfz_2,\by)\le L \big\}\prod_{i=1}^2 \one \{ x+r_ny_i\in Q_\ell \}  \\
&\qquad \qquad \qquad \qquad \times \Big(1-e^{-n\text{vol} \big(\U(x, x+r_ny_1, x + r_n y_2) \cap Q_\ell \big) }  \Big) \dif \by \dif x  \\
&\le \f1{6}\int_{(\R^2)^2} \one \big\{ \diam(\bfz_2,\by)\le L \big\} \Big(1-e^{-nr_n^2 \mathcal R(\bfz_2,\by)^2 \pi  }  \Big) \dif \by \to 0.
\end{align*}

Second, we also  need to modify the proof of  \eqref{e:UI1}. 
Specifically, we need to show that  
\begin{equation*}
\limsup_{n\to\infty}  \E \Big[ e^{a\sum_{\Y\subset \Pnr_{Q_1}} s_n( \Y, \Pnr_{Q_1} )   }  \Big]  <\infty. 
\end{equation*}
For this purpose, we again consider the diluted families of cubes $G_1,G_2,\dots$ as in \eqref{e:diluted.cubes} (with $t\equiv 1$). Then, as an analog of \eqref{e:equiv.spatial.indep}, our task is reduced to showing that for every $a>0$, 
\begin{equation}  \label{e:UI.critical.point}
\limsup_{n\to\infty} \bigg( \E \Big[ e^{a\sum_{\Y\subset \Pnr_{Q_1}} s_n( \Y, \Pnr_{Q_1} ) \one \{ \ell(\Y)\in J_1 \}  }  \Big]  \bigg)^{b_n}<\infty,  
\end{equation}
where $J_1=[0,r_n/\sqrt{2}]^2$ and $b_n=\rho_{3,n}^{-1}/(4Lr_n)^2$. Here, a key insight is that if $s_n(\Y, \Pnr_{Q_1}) = 1$ for some $\Y\subset \Pnr_{Q_1}$, then $\Y$ defines a Voronoi point in the plane. As the number of Voronoi points on a set of $m$ vertices in $\R^2$ is of at most $\mathcal O(m)$ (see \cite[Lemma 1]{lee:schachter:1980}), we can deduce that 
\begin{align}  
&\sum_{\Y\subset \Pnr_{Q_1}} s_n(\Y, \Pnr_{Q_1})\, \one \big\{ \ell(\Y)\in J_1 \big\} \label{e:Voronoi.points} \\
&\qquad \le C^* \big|  \Pn\cap \text{Tube}(J_1,r_nL) \big|\, \one \Big\{ \big|  \Pn\cap \text{Tube}(J_1,r_nL) \big|\ge3 \Big\}, \notag 
\end{align}
where 
$
\text{Tube}(J_1,r_nL)  := \big\{ x\in \R^2: \inf_{y\in J_1}\| x-y\| \le r_nL \big\}. 
$
By \eqref{e:Voronoi.points}, one can see that 
\begin{align*}
&\E\Big[ e^{a\sum_{\Y\subset \Pnr_{Q_1}} s_n(\Y, \Pnr_{Q_1})\, \one \{ \ell(\Y)\in J_1 \} } \Big] \\
&\le \E\Big[ e^{aC^*|  \Pn\cap \text{Tube}(J_1,r_nL) |\, \one \{ |  \Pn\cap \text{Tube}(J_1,r_nL) |\ge3 \}} \Big] \\
&=\sum_{\ell=3}^\infty e^{aC^*\ell} \P \Big( \big| \Pn \cap \text{Tube} (J_1,r_nL) \big|=\ell \Big)=\mathcal O\big( (nr_n^2)^3 \big). 
\end{align*}
This concludes the proof of \eqref{e:UI.critical.point}.
\end{proof}

\begin{acks}[Acknowledgments]
The authors are very grateful for useful comments received from an anonymous referee
and an anonymous Associate Editor. The referee proposed interesting topics for further research, while helping the authors to introduce
a number of improvements to the paper.
\end{acks}

\begin{funding}
TO's research was supported by the NSF grant DMS-1811428 and the AFOSR grant FA9550-22-0238. CH would like to acknowledge the financial support of the CogniGron research center and the Ubbo Emmius Funds (Univ. of Groningen).
\end{funding}
\bibliographystyle{imsart-number} 
\bibliography{LDP-Betti}

\end{document}